\documentclass[11pt,a4paper]{article}

\usepackage{CSPpreamble}

\newtheorem{example}{Example}[section]

\numberwithin{equation}{section}

\begin{document}

\title{A Self-Consistent Field Solution \\
for Robust Common Spatial Pattern Analysis}
\author{
Dong Min Roh\thanks{Department of Mathematics, University of California, Davis, CA 95616, USA. \href{droh@ucdavis.edu}{droh@ucdavis.edu}}
~and~
Zhaojun Bai\thanks{Department of Computer Science and Department of Mathematics, University of California, Davis, CA 95616, USA. \href{zbai@ucdavis.edu}{zbai@ucdavis.edu}}
}
\date{}
\maketitle

\vspace{-1em}

\begin{abstract}

The common spatial pattern analysis (CSP) is a widely used signal processing technique in brain-computer interface (BCI) systems to increase the signal-to-noise ratio in electroencephalogram (EEG) recordings.
Despite its popularity, the CSP's performance is often hindered by the nonstationarity and artifacts in EEG signals.
The minmax CSP improves the robustness of the CSP by using data-driven covariance matrices to accommodate the uncertainties. 
We show that by utilizing the optimality conditions, the minmax CSP can be recast as an eigenvector-dependent nonlinear eigenvalue problem (NEPv).
We introduce a self-consistent field (SCF) iteration with line search that solves the NEPv of the minmax CSP.
Local quadratic convergence of the SCF for solving the NEPv is illustrated using synthetic datasets. 
More importantly, experiments with real-world EEG datasets show the improved motor imagery classification rates and shorter running time of the proposed SCF-based solver compared to the existing algorithm for the minmax CSP.
\end{abstract}



\section{Introduction}

\paragraph{BCI, EEG, and ERD.}
Brain-computer interface (BCI) is a system that translates the participating subject's intent into control actions, such as control of computer applications or prosthetic devices \cite{dornhege2007toward,wolpaw2002brain,lemm2011introduction}.
BCI studies have utilized a range of approaches for capturing the brain signals pertaining to the subject's intent, from invasive (such as brain implants) and partially invasive (inside the skull but outside the brain) to non-invasive (such as head-worn electrodes). 
Among non-invasive approaches, electroencephalogram (EEG) based BCIs are the most widely studied due to their ease of setup and non-invasive nature. 

In EEG based BCIs, the brain signals are often captured through the electrodes (small metal discs) attached to the subject's scalp while they perform motor imagery tasks. 
Two phases, the calibration phase and the feedback phase, typically comprise an EEG-BCI session \cite{blankertz2007optimizing}. 
During the calibration phase, the preprocessed and spatially filtered EEG signals are used to train a classifier.
In the feedback phase, the trained classifier interprets the subject's EEG signals and translates them into control actions, such as cursor control or button selection, often occurring in real-time.

Translation of a motor imagery into commands is based on the phenomenon known as event-related desynchronization (ERD) \cite{pfurtscheller1999event, blankertz2007optimizing}, a decrease in the rhythmic activity over the sensorimotor cortex as a motor imagery takes place.
Specific locations in the sensorimotor cortex are related to corresponding parts of the body. 
For example, the left and right hands correspond to the right and left motor cortex, respectively. 
Consequently, the cortical area where an ERD occurs indicates the  specific motor task being performed.

\paragraph{CSP.} 
Spatial filtering is a crucial step to uncover the discriminatory modulation of an ERD \cite{blankertz2007optimizing}. 
One popular spatial filtering technique in EEG-based BCIs is the common spatial pattern analysis (CSP) \cite{koles1991quantitative,blankertz2007optimizing,lotte2010regularizing,reuderink2008robustness}. 
The CSP is a data-driven approach that computes subject-specific spatial filters to increase the signal-to-noise ratio in EEG signals and uncover the discrimination between two motor imagery conditions \cite{samek2012stationary}. 
Mathematically, the principal spatial filters that optimally discriminate the motor imagery conditions can be expressed as solutions to generalized Rayleigh quotient optimizations.

\paragraph{Robust CSP and existing work.}
However, the CSP is often subject to poor performance due to the nonstationary nature of EEG signals \cite{blankertz2007optimizing,samek2012stationary}.
The nonstationarity of a signal is influenced by various factors such as artifacts resulting from non-task-related muscle movements and eye movements, as well as between-session errors such as small differences in electrode positions and a gradual increase in the subject's tiredness. 
Consequently, some of the EEG signals can be noisy and act as outliers.
This leads to poor average covariance matrices that do not accurately capture the underlying variances of the signals.
As a result, the CSP, which computes its principal spatial filters using the average covariance matrices in the Rayleigh quotient, can be highly sensitive to noise and prone to overfitting \cite{reuderink2008robustness,yong2008robust}.

There are many variants of the CSP that seek to improve its robustness. 
One approach directly tackles the nonstationarities, which includes regularizing the CSP to ensure its spatial filters remain invariant to nonstationarities in the signals \cite{blankertz2007invariant, samek2012stationary} or by first removing the nonstationary contributions before applying the CSP \cite{samek2012brain}. 
Another approach includes the use of divergence functions from information geometry, such as the Kullback–Leibler divergence \cite{arvaneh2013optimizing} and the beta divergence \cite{samek2013robust}, to mitigate the influence of outliers.
Reformulations of the CSP based on norms other than the L2-norm also seek to suppress the effect of outliers.
These norms include the L1-norm \cite{wang2011l1, li2016robust} and the L21-norm \cite{gu2021common}.
Another approach aims to enhance the robustness of the covariance matrices. 
This includes the dynamic time warping approach \cite{azab2019robust}, the application of the minimum covariance determinant estimator \cite{yong2008robust}, and the minmax CSP \cite{kawanabe2009robust, kawanabe2014robust}.

\paragraph{Contributions.} In this paper, we focus on the minmax CSP,  
a data-driven approach that utilizes the sets reflecting the regions 
of variability of covariance 
matrices~\cite{kawanabe2009robust,kawanabe2014robust}, 
The computation of spatial filters of the minmax CSP corresponds to solving the nonlinear Rayleigh quotient optimizations (OptNRQ). 
In \cite{kawanabe2009robust,kawanabe2014robust}, the OptNRQ is tackled by iteratively solving the linear Rayleigh quotient optimizations. 
However, this method is susceptible to convergence failure, as demonstrated in Section~\ref{subsec:conv}. Moreover, even upon convergence, a more fundamental issue persists: the possibility of obtaining a suboptimal solution.
We elucidate that this is due to a so-called
eigenvalue ordering issue associated with its underlying 
eigenvector-dependent nonlinear eigenvalue problem (NEPv) 
to be discussed in Section~\ref{subsec:fp}.

This paper builds upon an approach for solving a general OptNRQ proposed in \cite{bai2018robust}.
In \cite{bai2018robust}, authors studied the robust Rayleigh quotient optimization problem where the data matrices of the Rayleigh quotient are subject to uncertainties. 
They proposed to solve such a problem by exploiting its characterization as a NEPv. The application to the minmax CSP was outlined in their work. 
Our current work is a refinement and extension of the NEPv approach to solving the minmax CSP. 
(i) We provide a comprehensive theoretical background and step-by-step derivation 
of the NEPv formulation for the minmax CSP. 
(ii) We introduce an algorithm, called the OptNRQ-nepv, a self-consistent field (SCF) 
iteration with line search to solve the minmax CSP. 
(iii) Through synthetic datasets,  we illustrate 
the  ``eigenvalue-order issue'' to the convergence of algorithms and 
demonstrate the local quadratic convergence of 
the proposed OptNRQ-nepv algorithm. 
(iv) More importantly, the OptNRQ-nepv is applied to real-world EEG datasets. 
Numerical results show improved motor imagery classification rates 
of the OptNRQ-nepv to existing algorithms. 

\paragraph{Paper organization.}
The rest of this paper is organized as follows: 
In Section~\ref{sec:CSP}, we give an overview of the CSP, including the data processing steps and the Rayleigh quotient optimizations for computing the spatial filters. 
Moreover, the minmax CSP and its resulting nonlinear Rayleigh quotient optimizations are summarized.
In Section~\ref{sec:NEPv}, we provide a detailed derivation and explicit formulation of the NEPv arising from the minmax CSP, followed by an example to illustrate the eigenvalue ordering issue. Algorithms for solving the minmax CSP are presented in Section~\ref{sec:alg}.
Section~\ref{sec:num} consists of numerical experiments on synthetic and real-world EEG datasets, which include results on convergence behaviors, running time, and classification performance of algorithms. 
Finally, the paper concludes with concluding remarks in Section~\ref{sec:con}. 


\section{CSP and minmax CSP} \label{sec:CSP}

\subsection{Data preprocessing and data covariance matrices}

An EEG-BCI session involves collecting trials, which are time series of EEG signals recorded from electrode channels placed on the subject's scalp during motor imagery tasks. 
Preprocessing is an important step to improve the quality of the raw EEG signals in the trials, which leads to better performance of the CSP \cite{blankertz2007optimizing,dornhege200713,muller1999designing}.
The preprocessing steps include channel selection, bandpass filtering, and time interval selection of the trials. 
Channel selection involves choosing the electrode channels around the cortical area corresponding to the motor imagery actions. 
Bandpass filtering is used to filter out low and high frequencies that may originate from non-EEG related artifacts. 
Time interval selection is used to choose the time interval in which the motor imagery takes place. 

We use $Y_c^{(i)}$ to denote the $i^{\mbox{th}}$ trial for the motor imagery condition $c$.
Each trial $Y_c^{(i)}$ is a $n\times t$ matrix, whose row corresponds to an EEG signal sampled at an electrode channel with $t$ many sampled time points. 
We assume that the trials $Y_c^{(i)}$
are scaled and centered. Otherwise,
\[Y_c^{(i)}\longleftarrow\frac{1}{\sqrt{t-1}}Y_c^{(i)}(I_t-\frac{1}{t}\mathbbm{1}_t\mathbbm{1}_t^T).\]
Here, $I_t$ denotes the $t\times t$ identity matrix and $\mathbbm{1}_t$ denotes the $t$-dimensional vector of all ones. 
The data covariance matrix for the trial $Y_c^{(i)}$ is defined as
\begin{equation}\label{eq:data_cov}
\Sigma_c^{(i)}:=Y_c^{(i)}{Y_c^{(i)}}^T \in \R^{n\times n}.
\end{equation}
The average of the data covariance matrices is given by
\begin{equation}\label{eq:cov}
\widehat{\Sigma}_c
=\frac{1}{N_c}\sum_{i=1}^{N_c}\Sigma_c^{(i)} \in\R^{n\times n},
\end{equation}
where $N_c$ is the total number of trials of condition $c$.

\begin{assumption} \label{assume:sigma_data}
We assume that the average covariance matrix $\widehat{\Sigma}_c$~\eqref{eq:cov} is positive definite, i.e., 
\begin{equation}\label{eq:Avg_Sig_pos}
    \widehat{\Sigma}_c\succ0.
\end{equation}
\end{assumption}


\subsection{CSP}

\paragraph{Spatial filters.} 
In this paper, we consider the binary CSP analysis where two conditions 
of interest, such as left hand and right hand motor imageries, 
are represented by $c\in\{-,+\}$. 
With respect to $c$, we use $\bar{c}$ to denote its complement condition. That is, if $c=-$ then $\bar{c}=+$, and if $c=+$ then $\bar{c}=-$.
For the discussion on the CSP involving more than two conditions, we refer 
the readers to \cite{blankertz2007optimizing} and references therewithin.

The principle of the CSP is to find spatial filters such that the variance of the spatially filtered signals under one condition is minimized, while that of the other condition is maximized. 
Representing the $n$ spatial filters as the columns of $X\in\R^{n\times n}$, 
the following properties of the CSP 
hold \cite{koles1991quantitative,blankertz2007optimizing}:

\begin{enumerate}[label=(\roman*)]
    \item The average covariance matrices $\widehat{\Sigma}_c$ and 
$\widehat{\Sigma}_{\bar{c}}$ are simultaneously diagonalized:
    \begin{equation}\label{eq:simul_diag}
    \begin{split}
        & X^T\widehat{\Sigma}_cX=\Lambda_c, \\
        & X^T\widehat{\Sigma}_{\bar{c}}X=\Lambda_{\bar{c}},
    \end{split}
    \end{equation}
    where $\Lambda_c$ and $\Lambda_{\bar{c}}$ are diagonal matrices.

    \item The spatial filters $X$ are scaled such that 
    \begin{equation}\label{eq:identity_con}
        X^T(\widehat{\Sigma}_c + \widehat{\Sigma}_{\bar{c}})X=
        \Lambda_c+\Lambda_{\bar{c}}=I_n,
    \end{equation}
    where $I_n$ denotes the $n \times n$ identity matrix.
\end{enumerate}
Let $Z_c^{(i)}$ denote the spatially filtered trial of $Y_c^{(i)}$, i.e., 
\begin{equation}\label{eq:filt_sig}
    Z_c^{(i)}:=X^T{Y_c^{(i)}}.
\end{equation}
Then
\begin{enumerate}[label=(\roman*)]

    \item By the property~\eqref{eq:simul_diag}, the average covariance matrix of $Z_c^{(i)}$ is diagonalized:
    \begin{equation}\label{eq:filt_diag}
        \frac{1}{N_c}\sum_{i=1}^{N_c}Z_c^{(i)}{Z_c^{(i)}}^T=\frac{1}{N_c}\sum_{i=1}^{N_c}X^TY_c^{(i)}{Y_c^{(i)}}^TX=X^T\Big(\frac{1}{N_c}\sum_{i=1}^{N_c}\Sigma_c^{(i)}\Big)X=X^T\widehat{\Sigma}_cX=\Lambda_c.
    \end{equation}
    This indicates that the spatially filtered signals between electrode channels are uncorrelated.
    
    \item The diagonals of $\Lambda_c$~\eqref{eq:filt_diag} correspond to the average variances of the filtered signals in condition $c$. 
    Since both $\widehat{\Sigma}_c$ and $\widehat{\Sigma}_{\bar{c}}$ are positive definite by Assumption~\ref{assume:sigma_data}, the property~\eqref{eq:identity_con} indicates that the variance of either condition must lie between 0 and 1, and that the sum of the variances of the two conditions must be 1. 
    This indicates that a small variance in one condition implies a large variance in other condition, and vice versa.
\end{enumerate}

\paragraph{Connection of spatial filters and generalized eigenvalue problems.} 
The properties \eqref{eq:simul_diag} and \eqref{eq:identity_con} together imply 
that the spatial filters $X$ correspond to the eigenvectors of 
the generalized eigenvalue problem
\begin{equation}\label{eq:gen_eig}
\begin{cases}
\widehat{\Sigma}_c X = (\widehat{\Sigma}_c + \widehat{\Sigma}_{\bar{c}})X\Lambda_c \\
X^T (\widehat{\Sigma}_c+\widehat{\Sigma}_{\bar{c}}) X=I_n.
\end{cases}
\end{equation}

\paragraph{Principal spatial filters and Rayleigh quotient optimizations.}
In practice, only a small subset of the spatial filters, rather than 
all $n$ of them, are necessary \cite{kawanabe2014robust,blankertz2007optimizing}. 
In this paper, we consider the computations of the 
\textbf{principal spatial filters} $x_c$, namely, 
$x_c$ is the eigenvector of the matrix pair 
$(\widehat{\Sigma}_c,\widehat{\Sigma}_c + \widehat{\Sigma}_{\bar{c}})$ 
with respect to the smallest eigenvalue.
By the Courant-Fischer variational principle \cite{parlett1998symmetric}, 
$x_c$ is the solution of the following generalized Rayleigh quotient optimization 
(RQopt)\footnote{Strictly speaking, the solution of the Rayleigh quotient 
optimization~\eqref{eq:gRq} may not satisfy the normalization constraint 
$x^T(\widehat{\Sigma}_c+\widehat{\Sigma}_{\bar{c}})x=1$. But, in this case, 
we can always normalize the obtained solution by applying the transformation 
$x\longleftarrow\frac{x}{\sqrt{x^T(\widehat{\Sigma}_c+\widehat{\Sigma}_{\bar{c}})x}}$ 
to obtain $x_c$. Hence, we can regard the generalized Rayleigh quotient 
optimization~\eqref{eq:gRq} as the optimization of interest for computing 
the principal spatial filter $x_c$.}
\begin{equation}\label{eq:gRq}
\min_{x\neq0} \frac{x^T\widehat{\Sigma}_c x}
{x^T(\widehat{\Sigma}_c + \widehat{\Sigma}_{\bar{c}})x}.
\end{equation}

\subsection{Minmax CSP}\label{sec:rCSP}
The minmax CSP proposed in \cite{kawanabe2009robust,kawanabe2014robust} improves CSP's robustness to the nonstationarity and artifacts present in the EEG signals. 
Instead of being limited to fixed covariance matrices 
as in \eqref{eq:gRq}, the minmax CSP considers the sets $\mathcal{S}_c$ 
containing candidate covariance matrices. 
In this approach, the principal spatial filters $x_c$
are computed by solving the minmax optimizations\footnote{We point out that 
in \cite{kawanabe2009robust,kawanabe2014robust}, the optimizations are 
stated as equivalent maxmin optimizations. For instance, $\max_{x\neq0}
\min_{\substack{\Sigma_c\in\mathcal{S}_c\\
\Sigma_{\bar{c}}\in\mathcal{S}_{\bar{c}}}}  
(x^T\Sigma_{\bar{c}}x)/(x^T(\Sigma_c+\Sigma_{\bar{c}})x)$ is equivalent to 
our minmax optimization~\eqref{eq:cspminmax-} using 
$(x^T\Sigma_{\bar{c}}x)/(x^T(\Sigma_c+\Sigma_{\bar{c}})x)
=1-(x^T\Sigma_cx)/(x^T(\Sigma_c+\Sigma_{\bar{c}})x)$.
}
\begin{equation}\label{eq:cspminmax-}
\min_{x\neq0}
\max_{
       \substack{\Sigma_c\in\mathcal{S}_c \\ 
                 \Sigma_{\bar{c}}\in\mathcal{S}_{\bar{c}}
                }  
     } 
    \frac{x^T\Sigma_cx}{x^T(\Sigma_c+\Sigma_{\bar{c}})x}.
\end{equation}
The sets $\mathcal{S}_c$ are called the {\bf tolerance sets} and are constructed to reflect the regions of variability of the covariance matrices.
The inner maximization in \eqref{eq:cspminmax-} corresponds to seeking the worst-case generalized Rayleigh quotient within all possible covariance matrices in the tolerance sets. 
Optimizing for the worst-case behavior is a popular paradigm in optimization to obtain solutions that are more robust to noise and overfitting \cite{ben2009robust}.

\paragraph{Tolerance sets of the minmax CSP.} 
Two approaches for constructing the tolerance sets $\mathcal{S}_c$ are proposed in \cite{kawanabe2009robust,kawanabe2014robust}. 
\begin{enumerate}[label=(\roman*)]
    \item The first approach defines $\mathcal{S}_c$ as the set of covariance matrices that are positive definite and bounded around the average covariance matrices by a weighted Frobenius norm. 
    This approach leads to generalized Rayleigh quotient optimizations, similar to the CSP (see Appendix~\ref{app:Frob} for details).
    While this approach is shown to improve the CSP in classification, 
the corresponding tolerance sets do not adequately capture the 
true variability of the covariance matrices over time.
    Additionally, the suitable selection of positive definite matrices for the weighted Frobenius norm remains an open problem. 
    These limitations hinder the practical application of this approach.
    
    \item The second approach is data-driven and defines $\mathcal{S}_c$ as the set of covariance matrices described by the interpolation matrices derived from the data covariance matrices $\{\Sigma_c^{(i)}\}_{i=1}^{N_c}$~\eqref{eq:data_cov}. 
    The tolerance sets in this approach can effectively capture the variability present in the data covariance matrices and eliminate the need for additional assumptions or prior information.
    This data-driven approach leads to nonlinear Rayleigh quotient optimizations. We will focus on this approach.
    
\end{enumerate}

\paragraph{Data-driven tolerance sets.}
Given the average covariance matrix $\widehat\Sigma_{c}$~\eqref{eq:cov}, the symmetric interpolation matrices $V_c^{(i)}\in\R^{n\times n}$ for $i= 1,\ldots,m$, and the weight matrix $W_c\in\R^{m\times m}$ with positive weights $w_c^{(i)}$,
\begin{equation}\label{eq:Wc}
    W_c := \mbox{diag}\big(w_c^{(1)}, w_c^{(2)}, \ldots, w_c^{(m)}\big),
\end{equation}
the tolerance sets of the data-driven approach are constructed as 
\begin{align}\label{eq:cspset}
\mathcal{S}_{c } := \left\{ 
\Sigma_{c}(\alpha_c) = 
\widehat\Sigma_{c} +\sum_{i=1}^{m}\alpha_{c}^{(i)} V_{c}^{(i)}
\,\,\bigg|\,\,
\|\alpha_c\|_{W_c^{-1}}\leq\delta_c
\right\},
\end{align}
where $\delta_c>0$ denotes the tolerance set radius.
$\alpha_c = [\alpha_{c}^{(1)}, \alpha_{c}^{(2)}, \ldots, \alpha_{c}^{(m)}]^T$, and $\|\alpha_c\|_{W_c^{-1}}$ is a weighted vector norm, defined with matrix $W_c^{-1}$ as\footnote{In \cite{kawanabe2009robust,kawanabe2014robust}, the norm~\eqref{eq:weighted_norm} is referred to as the \textit{PCA-based norm}, reflecting the derivation of the weights $\{w_c^{(i)}\}_{i=1}^m$ from the PCA applied to the data covariance matrices. However, we recognize the norm~\eqref{eq:weighted_norm} as a weighted vector norm defined for general values of $w_c^{(i)}$.}
\begin{equation}\label{eq:weighted_norm} \|\alpha_c\|_{W_c^{-1}}:=\sqrt{\alpha_c^TW_c^{-1}\alpha_c}.
\end{equation}
The symmetric interpolation matrices $\{V_c^{(i)}\}_{i=1}^{m}$ and the weight matrix $W_c$ are parameters that are obtained from the data covariance matrices 
$\{\Sigma_c^{(i)}\}_{i=1}^{N_c}$ to reflect their variability.

\paragraph{Why $\Sigma_c(\alpha_{c})$ of $\mathcal{S}_c$~\eqref{eq:cspset} increases robustness.}
The average covariance matrix $\widehat{\Sigma}_c$ tends to overestimate the largest eigenvalue and underestimate the smallest eigenvalue \cite{ledoit2004well}. 
For the CSP, which relies on the average covariance matrices, this effect can significantly influence the variance of two conditions.
On the other hand, the covariance matrices $\Sigma_c(\alpha_{c})$ of the tolerance set $\mathcal{S}_c$~\eqref{eq:cspset} reduce the relative impact of under- and over-estimations by the addition of the term $\Sigma_{i=1}^m\alpha_c^{(i)}V_c^{(i)}$ to $\widehat{\Sigma}_c$. 
As a result, $\Sigma_c(\alpha_{c})$ can robustly represent the variance of its corresponding condition.
Another factor contributing to the robustness of $\Sigma_c(\alpha_c)$ is the adequate size of the tolerance set $\mathcal{S}_c$.
If the radius $\delta_c$ of $\mathcal{S}_c$ is zero, then $\Sigma_c(\alpha_c)$ is equal to $\widehat{\Sigma}_c$, the covariance matrix used in the standard CSP, and no variability is taken into account.
If $\delta_c$ is too large, then $\Sigma_c(\alpha_c)$ is easily influenced by the outliers.
On the other hand, by selecting an appropriate $\delta_c$, $\Sigma_c(\alpha_{c})$ can capture the intrinsic variability of the data while ignoring outliers that may dominate the variability.

\paragraph{Computing $V^{(i)}_c$ and $W_c$.}
One way to obtain the interpolation matrices $V^{(i)}_c$ and the weight matrix 
$W_c$ is by performing the principal component analysis (PCA) on the 
data covariance matrices $\{\Sigma_c^{(i)}\}_{i=1}^{N_c}$.
The following steps outline this process:
\begin{enumerate}
    \item Vectorize each data covariance matrix $\Sigma_c^{(i)}$ by stacking its columns into a $n^2$-dimensional vector, i.e., obtain $\mbox{vec}(\Sigma_c^{(i)})\in\R^{n^2}$.
    \item Compute the covariance matrix $\Gamma_c\in\R^{n^2\times n^2}$ of the vectorized covariance matrices $\{\mbox{vec}(\Sigma_c^{(i)})\}_{i=1}^{N_c}$.
    \item Compute the $m$ largest eigenvalues and corresponding eigenvectors (principal components) of $\Gamma_c$ as $\{w_c^{(i)}\}_{i=1}^{m}$ and $\{\nu_c^{(i)}\}_{i=1}^{m}$, respectively.
    \item Transform the eigenvectors $\{\nu_c^{(i)}\}_{i=1}^{m}$ back into $n\times n$ matrices, i.e., obtain $\mbox{mat}(\nu_c^{(i)}) \in \R^{n \times n}$, then symmetrizing afterwards to obtain the interpolation matrices $\{V_c^{(i)}\}_{i=1}^{m}$.
    \item 
    Form the weight matrix $W_c := \mbox{diag}\big(w_c^{(1)}, w_c^{(2)}, \ldots, w_c^{(m)}\big)$.
\end{enumerate}
The interpolation matrices $V^{(i)}_c$ and the weights $w_c^{(i)}$ obtained in this process correspond to the principal components and variances of the data covariance matrices $\{\Sigma_c^{(i)}\}_{i=1}^{N_c}$, respectively.
Consequently, the tolerance set $\mathcal{S}_c$~\eqref{eq:cspset} has a geometric interpretation as an ellipsoid of covariance matrices that is centered at the average covariance matrix $\widehat{\Sigma}_c$ and spanned by the principal components $V^{(i)}_c$.
Moreover, a larger variance $w_c^{(i)}$ allows a larger value of $\alpha_c^{(i)}$ according to the norm~\eqref{eq:weighted_norm}.
This implies that more significant variations are permitted in the directions indicated by the leading principal components.
Therefore, the tolerance set $\mathcal{S}_c$ constructed by this PCA process can effectively capture the variability of the data covariance matrices.

We point out that the covariance matrix $\Gamma_c$ may only be positive semi-definite if the number of data covariance matrices is insufficient. 
In such cases, some of the eigenvalues $w_c^{(i)}$ may be zero, causing the norm definition~\eqref{eq:weighted_norm} to be ill-defined. 
To address this issue, we make sure to choose $m$ such that all $\{w_c^{(i)}\}_{i=1}^m$ are positive.

\paragraph{Minmax CSP and nonlinear Rayleigh quotient optimization.}\label{subsec:CSP_NRQopt}
We first note that a covariance matrix $\Sigma_c(\alpha_c)$ of the tolerance 
set $\mathcal{S}_c$~\eqref{eq:cspset} can be characterized as a function on 
the set
\begin{equation}\label{eq:omg_set}
\Omega_c:=\left\{\alpha_c\in\R^m \,\big|\,
\|\alpha_c\|_{W_c^{-1}}\leq\delta_c\right\}.
\end{equation}
Consequently, we can express the optimization \eqref{eq:cspminmax-} 
of the minmax CSP as
\begin{align}
\min_{x\neq0}
\max_{ \substack{
        \alpha_c\in \Omega_{c}\\
        \alpha_{\bar{c}}\in \Omega_{\bar{c}}
         }
     }  
    \frac{x^T\Sigma_c(\alpha_c)x}
         {x^T(\Sigma_c(\alpha_c)+\Sigma_{\bar{c}}(\alpha_{\bar{c}}))x}. 
\label{eq:cspminmax-2}
\end{align} 

\begin{assumption}\label{assume:sigma_alpha}
We assume that each covariance matrix in the tolerance set~\eqref{eq:cspset} is positive definite, i.e., 
\begin{equation}\label{eq:cov_pos}
\Sigma_c(\alpha_c)\succ0 
\quad \mbox{for all $\alpha_c\in\Omega_c$}.
\end{equation}
\end{assumption}
Assumption~\ref{assume:sigma_alpha} holds in practice for a small tolerance set radius $\delta_c$. 
If $\Sigma_c(\alpha_c)$ fails the positive definiteness property, we may set its negative eigenvalues to a small value, add a small perturbation term $\epsilon I$, or perform the modified Cholesky algorithm (as described in \cite{cheng1998modified}) to obtain a positive definite matrix. 


\paragraph{Nonlinear Rayleigh quotient optimizations.}
In the following, we show that the minmax CSP \eqref{eq:cspminmax-2} 
can be formulated as a nonlinear Rayleigh quotient optimization (OptNRQ) \cite{kawanabe2009robust,kawanabe2014robust,bai2018robust}. 

\begin{theorem}\label{thm:rbst_csp}
Under Assumption~\ref{assume:sigma_alpha}, the minmax CSP \eqref{eq:cspminmax-2} is equivalent to the following 
OptNRQ: 
\begin{equation}\label{eq:cspnrq}
\min_{x\neq0}\bigg\{q_c(x):=\frac{x^T\Sigma_c(x)x}{x^T(\Sigma_c(x)+\Sigma_{\bar{c}}(x))x}\bigg\}
\end{equation}
where
\begin{equation}\label{eq:Sigma_alpha_x}
\Sigma_{c}(x) = \widehat\Sigma_{c}
+\sum_{i=1}^m\alpha_{c}^{(i)}(x) V_{c}^{(i)} \succ 0.
\end{equation}
$\alpha_c^{(i)}(x)$ is a component of $\alpha_c(x):\R^n\to\R^m$ defined as
\begin{equation}\label{eq:alph_x}
\alpha_c(x) := \frac{\delta_c}{\|v_c(x)\|_{W_c}}W_c v_c(x),
\end{equation}
where $W_c$ is the weight matrix defined in~\eqref{eq:Wc} and
$v_c(x):\R^n\to\R^m$ is a vector-valued function defined as
\begin{equation}\label{eq:vcx}
    v_c(x) := 
    \begin{bmatrix}
        x^TV_c^{(1)}x, & x^TV_c^{(2)}x, & \ldots, & x^TV_c^{(m)}x
    \end{bmatrix}^T.
\end{equation}
$\Sigma_{\bar{c}}(x)$ is defined identically as 
$\Sigma_c(x)$~\eqref{eq:Sigma_alpha_x} with $c$ replaced with $\bar{c}$, 
except that
\begin{equation}\label{eq:alph_x2}
\alpha_{\bar{c}}(x) := -\frac{\delta_{\bar{c}}}{\|v_{\bar{c}}(x)\|_{W_{\bar{c}}}}W_{\bar{c}} v_{\bar{c}}(x)
\end{equation}
differs from $\alpha_c(x)$~\eqref{eq:alph_x} by a sign.
\end{theorem}
\begin{proof}
By Assumption~\ref{assume:sigma_alpha}, the inner maximization of the minmax CSP \eqref{eq:cspminmax-2} is equivalent to solving two separate optimizations:
\begin{equation}\label{eq:realinnmax}
    \max_{\substack{\alpha_c\in \Omega_{c}\\
        \alpha_{\bar{c}}\in \Omega_{\bar{c}}}}  
    \frac{x^T\Sigma_c(\alpha_c)x}{x^T(\Sigma_c(\alpha_c)+\Sigma_{\bar{c}}(\alpha_{\bar{c}}))x}
    =\frac{\max\limits_{\alpha_c\in \Omega_{c}}x^T\Sigma_c(\alpha_c)x}{\max\limits_{\alpha_c\in \Omega_{c}}x^T\Sigma_c(\alpha_c)x+\min\limits_{\alpha_{\bar{c}}\in \Omega_{\bar{c}}}x^T\Sigma_{\bar{c}}(\alpha_{\bar{c}})x}.
\end{equation}
We have
$$x^T\Sigma_c(\alpha_c)x=x^T\widehat{\Sigma}_cx+\sum_{i=1}^m\alpha_c^{(i)}(x^TV_c^{(i)}x)=x^T\widehat{\Sigma}_cx+\alpha_c^Tv_c(x)$$
with $v_c(x)$ as defined in~\eqref{eq:vcx}.
Moreover, 
\begin{align}
|\alpha_c^Tv_c(x)| &= 
|(W_c^{-1/2}\alpha_c)^T(W_c^{1/2}v_c(x))| \nonumber \\
&\leq \|W_c^{-1/2}\alpha_c\|_2\cdot\|W_c^{1/2}v_c(x)\|_2 
\label{eq:thm1_ineq1} \\
&\leq \delta_c \cdot \|v_c(x)\|_{W_c} \label{eq:thm1_ineq2}
\end{align}
where the inequality~\eqref{eq:thm1_ineq1} follows from the Cauchy-Schwarz inequality and the inequality~\eqref{eq:thm1_ineq2} follows from the constraint $\|\alpha_c\|_{W_c^{-1}}\leq\delta_c$ of the subspace $\Omega_c$~\eqref{eq:omg_set}.
With $\alpha_c$ defined as $\alpha_c(x)$ in~\eqref{eq:alph_x}, $|\alpha_c^Tv_c(x)|=\delta_c \cdot \|v_c(x)\|_{W_c}$.
Thus, $\alpha_c(x)$ is the solution to the term 
$\max\limits_{\alpha_c\in \Omega_{c}}x^T\Sigma_c(\alpha_c)x$.
By analogous arguments, $\alpha_{\bar{c}}$ defined as $\alpha_{\bar{c}}(x)$ in~\eqref{eq:alph_x2} is the solution to $\min\limits_{\alpha_{\bar{c}}\in \Omega_{\bar{c}}}x^T\Sigma_{\bar{c}}(\alpha_{\bar{c}})x$.
\end{proof}

\paragraph{Equivalent constrained optimizations.}
Due to the homogeneity of $\alpha_c(x)$~\eqref{eq:alph_x}, the covariance matrix 
$\Sigma_c(x)$~\eqref{eq:Sigma_alpha_x} is a homogeneous function, i.e., 
\begin{equation}\label{eq:Sigma_hom}
\Sigma_c(\gamma x)=\Sigma_c(x)
\end{equation}
for any nonzero $\gamma\in\R$.
The homogeneous property~\eqref{eq:Sigma_hom} also holds for condition $\bar{c}$.
This enables 
the OptNRQ~\eqref{eq:cspnrq} to be described by the following 
equivalent constrained optimization 
\begin{equation}\label{eq:cspnrq_new}
\begin{cases}
\min \quad  x^T\Sigma_c(x)x, \\ 
\mbox{s.t.} \quad x^T(\Sigma_c(x)+\Sigma_{\bar{c}}(x))x=1. 
\end{cases}
\end{equation}
Alternatively, the OptNRQ~\eqref{eq:cspnrq} can be described by the 
following equivalent optimization on the unit sphere $S^{n-1}:=\{x\in\R^n:x^Tx=1\}$:
\begin{equation}\label{eq:cspnrq_man}
\min_{x\in S^{n-1}}\frac{x^T\Sigma_c(x)x}{x^T(\Sigma_c(x)+\Sigma_{\bar{c}}(x))x}.
\end{equation}

While the constrained optimizations \eqref{eq:cspnrq_new} and \eqref{eq:cspnrq_man} are equivalent in the sense that an appropriate scaling of the solution of one optimization leads to the solution of the other, the algorithms for solving the two optimizations differ. 

In the next section, we start with the optimization~\eqref{eq:cspnrq_new} and show that it can be tackled by solving an associated eigenvector-dependent nonlinear eigenvalue problem (NEPv) using the self-consistent field (SCF) iteration.
The optimization~\eqref{eq:cspnrq_man} can be addressed using matrix manifold optimization methods (refer to \cite{absil2008optimization} for a variety of algorithms).  
In Section~\ref{subsec:conv}, we provide examples of convergence for the SCF iteration and the two matrix manifold optimization methods, the Riemannian conjugate gradient (RCG) algorithm and the Riemannian trust region (RTR) algorithm.
We demonstrate that the SCF iteration proves to be more efficient than the two manifold optimization algorithms.


\section{NEPv Formulations}\label{sec:NEPv}

In this section, we begin with a NEPv formulation of the OptNRQ~\eqref{eq:cspnrq} by utilizing the first order optimality conditions of the constrained optimization~\eqref{eq:cspnrq_new}. 
Unfortunately, this NEPv formulation encounters a so-called ``eigenvalue ordering issue''.\footnote{Given an eigenvalue $\lambda$ of a matrix, we use the terminology \textit{eigenvalue ordering} to define the position of $\lambda$ among all the eigenvalues $\lambda_1\leq\lambda_2\leq\cdots\leq\lambda_n$.} 
To resolve this issue, we derive an alternative NEPv by considering the second order optimality conditions of \eqref{eq:cspnrq_new}.

\subsection{Preliminaries}

In the following lemma, we present the gradients of functions $v_c(x),\|v_c(x)\|_{W_c}$ and $\alpha_c(x)$ with respect to $x$.

\begin{lemma}\label{lem:misc}
Let $v_c(x)\in\R^m$ and $\alpha_c(x)\in\R^m$ be as defined in \eqref{eq:vcx} and \eqref{eq:alph_x}. Then
    \begin{enumerate}[label=(\alph*)]
        \item $\nabla v_c(x)=
        2 \begin{bmatrix}
            V_c^{(1)}x, & V_c^{(2)}x, & \cdots, & V_c^{(m)}x
        \end{bmatrix}.$

        \item $\nabla \|v_c(x)\|_{W_c} = \frac{1}{\delta_c} \nabla v_c(x) \alpha_c(x)$.

        \item $\nabla \alpha_c(x) = -\frac{1}{\delta_c\|v_c(x)\|_{W_c}}\nabla v_c(x) \alpha_c(x)\alpha_c(x)^T + \frac{\delta_c}{\|v_c(x)\|_{W_c}} \nabla v_c(x)W_c$.
    \end{enumerate}
\end{lemma}
\begin{proof}
(a) For $i=1,2,\ldots,m$, the $i^{\mbox{th}}$ column of $\nabla v_c(x)\in\R^{n\times m}$ is $\nabla v_c^{(i)}(x)\in\R^n$, the gradient of the $i$-th component of $v_c(x)$. By Lemma~\ref{lem:der1}(b), $\nabla v_c^{(i)}(x)=\nabla(x^TV_c^{(i)}x)=2V_c^{(i)}x$, and the result follows.

(b) By definitions, \begin{align}
            \nabla\|v_c(x)\|_{W_c} &= \nabla\sqrt{v_c(x)^TW_cv_c(x)} = \frac{\nabla(v_c(x)^TW_cv_c(x))}{2\sqrt{v_c(x)^TW_cv_c(x)}} = \frac{2\nabla v_c(x)W_cv_c(x)}{2\sqrt{v_c(x)^TW_cv_c(x)}} \label{eq:vcx_wc_der1} \\
            &= \frac{1}{\delta_c} \nabla v_c(x)\alpha_c(x) \label{eq:vcx_wc_der2}
        \end{align}
        where we applied the chain rule and Lemma~\ref{lem:der2}(a) in~\eqref{eq:vcx_wc_der1}, and the definition of $\alpha_c(x)$~\eqref{eq:alph_x} in~\eqref{eq:vcx_wc_der2}.

(c) We apply Lemma~\ref{lem:der2}(c) with $a(x)=\frac{\delta_c}{\|v_c(x)\|_{W_c}}$ and $y(x)=W_cv_c(x)$ to obtain
        \begin{equation}\label{eq:alp_t_der2}
            \nabla\alpha_c(x) = \nabla\bigg(\frac{\delta_c}{\|v_c(x)\|_{W_c}}\bigg)\bigg(v_c(x)^TW_c\bigg) + \bigg(\frac{\delta_c}{\|v_c(x)\|_{W_c}}\bigg) \nabla\bigg(W_cv_c(x)\bigg).
        \end{equation}
        By the chain rule and applying the result $\nabla\|v_c(x)\|_{W_c}$ in Lemma~\ref{lem:misc}(b), we have
        \begin{equation}\label{eq:alp_t_der_first}
            \nabla\bigg(\frac{\delta_c}{\|v_c(x)\|_{W_c}}\bigg) = -\frac{1}{\|v_c(x)\|_{W_c}^2}\nabla v_c(x)\alpha_c(x).
        \end{equation}
        For $\nabla(W_cv_c(x))$, note that 
        $$W_cv_c(x)=\begin{bmatrix}
            w^{(1)}_c(x^TV_c^{(1)}x), &
            w^{(2)}_c(x^TV_c^{(2)}x), &
            \cdots, &
            w^{(m)}_c(x^TV_c^{(m)}x)
        \end{bmatrix}^T,$$
         so that its gradient is
        \begin{equation}\label{eq:alp_t_der_second}
            \nabla(W_cv_c(x))=\begin{bmatrix}
            2w_c^{(1)}V_c^{(1)}x & 2w_c^{(2)}V_c^{(2)}x & \cdots & 2w_c^{(m)}V_c^{(m)}x
        \end{bmatrix} = \nabla v_c(x)W_c
        \end{equation}
        by using $\nabla v_c(x)$ in Lemma~\ref{lem:misc}(a). 
        
        With results \eqref{eq:alp_t_der_first} and \eqref{eq:alp_t_der_second}, $\nabla\alpha_c(x)$~\eqref{eq:alp_t_der2} is
        \begin{align}
            \nabla\alpha_c(x) &= -\frac{1}{\|v_c(x)\|_{W_c}^2}\nabla v_c(x)\alpha_c(x)v_c(x)^TW_c + \frac{\delta_c}{\|v_c(x)\|_{W_c}}\nabla v_c(x)W_c \nonumber \\
            &= -\frac{1}{\delta_c\|v_c(x)\|_{W_c}}\nabla v_c(x)\alpha_c(x)\alpha_c(x)^T + \frac{\delta_c}{\|v_c(x)\|_{W_c}}\nabla v_c(x)W_c. \label{eq:alp_t_der_result}
        \end{align}
        where the equality \eqref{eq:alp_t_der_result} follows from the definition of $\alpha_c(x)$~\eqref{eq:alph_x}.
\end{proof}

We introduce the following quantity $s_c(x)$ for representing the quadratic form of the matrix $\Sigma_c(x)$:
\begin{equation}\label{eq:scx}
s_c(x):=x^T\Sigma_c(x)x.
\end{equation}

\begin{lemma}\label{lem:scx_der}
    Let $s_c(x)$ be the quadratic form as defined in~\eqref{eq:scx}. Then, its gradient is
    \begin{equation}\label{eq:scx_der}
        \nabla s_c(x)=2\Sigma_c(x)x.
    \end{equation}
\end{lemma}

\begin{proof}
    By the definition of $v_c(x)$~\eqref{eq:vcx}, we have 
\begin{align}\label{eq:sc_sim}
s_c(x)&=x^T\widehat{\Sigma}_cx + \sum_{i=1}^m\alpha_c^{(i)}(x)(x^TV_c^{(i)}x)=x^T\widehat{\Sigma}_cx + \alpha_c(x)^Tv_c(x) \nonumber \\
&= x^T\widehat{\Sigma}_cx + \delta_c\|v_c(x)\|_{W_c}
\end{align}
where the second term $\alpha_c(x)^Tv_c(x)$ was simplified to $\delta_c\|v_c(x)\|_{W_c}$. 
    
By Lemma~\ref{lem:der1}(b), the gradient of the first term $x^T\widehat{\Sigma}_cx$ of~\eqref{eq:sc_sim} is $2\widehat{\Sigma}_cx$. For the second term, we use the result $\nabla \|v_c(x)\|_{W_c} = \frac{1}{\delta_c} \nabla v_c(x)\alpha_c(x)$ in Lemma~\ref{lem:misc}(b).
    Therefore, by combining the gradients of the two terms in~\eqref{eq:sc_sim}, the gradient of $s_c(x)$ is
\begin{align}
 \nabla s_c(x) &= 2\widehat{\Sigma}_cx + \nabla v_c(x)\alpha_c(x) \label{eq:scx_der0} \\ 
 &= 2\widehat{\Sigma}_cx + 2\sum_{i=1}^m\alpha_c^{(i)}(x)V_c^{(i)}x \label{eq:scx_der1} \\
 &= 2\Sigma_c(x)x \label{eq:scx_der2}
 \end{align}
 where the equality \eqref{eq:scx_der1} follows from the result $\nabla v_c(x)$ in Lemma~\ref{lem:misc}(a) and the equality \eqref{eq:scx_der2} follows from the definition of $\Sigma_c(x)$~\eqref{eq:Sigma_alpha_x}.
\end{proof}


\begin{lemma}\label{lem:scx_der_der}
Let $s_c(x)$ be the quadratic form as defined in~\eqref{eq:scx}. Then, its Hessian matrix is
\begin{equation}\label{eq:scx_der_der}
\nabla^2s_c(x)=2(\Sigma_c(x)+\widetilde{\Sigma}_c(x)),
\end{equation}
where
\begin{equation}\label{eq:tilde_sig}
\widetilde{\Sigma}_c(x):=\frac{1}{2}\nabla\alpha_c(x)\nabla v_c(x)^T
\end{equation}
with $\nabla\alpha_c(x)$ and $\nabla v_c(x)$ as defined in Lemma~\ref{lem:misc}(c) and Lemma~\ref{lem:misc}(a), respectively.
\end{lemma}

\begin{proof}
    We compute the gradients of the two terms of $\nabla s_c(x)$ in~\eqref{eq:scx_der0} to compute the Hessian matrix $\nabla^2s_c(x)$. The gradient of the first term $2\widehat{\Sigma}_cx$ of \eqref{eq:scx_der0} is $2\widehat{\Sigma}_c$ by Lemma~\ref{lem:der1}(a). For the second term, we apply Lemma~\ref{lem:der2}(b) with $B(x)=\nabla v_c(x)$ and $y(x)=\alpha_c(x)$:
    $$\nabla(\nabla v_c(x)\alpha_c(x))=\sum_{i=1}^m(\alpha_c^{(i)}(x))(2V_c^{(i)}) + \nabla \alpha_c(x)\nabla v_c(x)^T.$$
    Therefore, by combining the gradients of the two terms in~\eqref{eq:scx_der0} together, we have
    \begin{align}
        \nabla^2 s_c(x) &= 2\widehat{\Sigma}_c + 2 \sum_{i=1}^m\alpha_c^{(i)}(x)V_c^{(i)} + \nabla \alpha_c(x)\nabla v_c(x)^T \nonumber \\
        &= 2\bigg(\widehat{\Sigma}_c + \sum_{i=1}^m\alpha_c^{(i)}(x)V_c^{(i)} + \frac{1}{2} \nabla \alpha_c(x)\nabla v_c(x)^T\bigg) \nonumber \\
        &= 2(\Sigma_c(x) + \widetilde{\Sigma}_c(x)) \label{eq:Hess_scx}
    \end{align}
where the equality~\eqref{eq:Hess_scx} follows from the definition of $\Sigma_c(x)$~\eqref{eq:Sigma_alpha_x} and $\widetilde{\Sigma}_c(x)$~\eqref{eq:tilde_sig}.  \end{proof}

\begin{remark}
    Due to the negative sign present in $\alpha_{\bar{c}}(x)$ for condition $\bar{c}$, the gradients $\nabla\|v_{\bar{c}}(x)\|_{W_{\bar{c}}}$ and $\nabla \alpha_{\bar{c}}(x)$ differ by a sign from the results in Lemma~\ref{lem:misc}(b) and Lemma~\ref{lem:misc}(c), respectively, when $c$ is replaced with $\bar{c}$.
    Regardless, the results in Lemma~\ref{lem:scx_der} and Lemma~\ref{lem:scx_der_der} hold for $\bar{c}$, i.e., $\nabla s_{\bar{c}}(x)=2\Sigma_{\bar{c}}(x)x$ and $\nabla^2s_{\bar{c}}(x)=2(\Sigma_{\bar{c}}(x)+\widetilde{\Sigma}_{\bar{c}}(x))$.
\end{remark}

\begin{lemma}\label{lem:Gcx}
Let 
\begin{equation}\label{eq:Gcx}
H_c(x) := \frac{1}{2} \nabla^2 s_c(x) = \Sigma_c(x)+\widetilde{\Sigma}_c(x).
\end{equation}
    Then, the following properties for $H_c(x)$ hold:
    \begin{enumerate}[label=(\alph*)]
        \item $H_c(x)$ is symmetric.
        \item $H_c(x)x=\Sigma_c(x)x$.
        \item Under Assumption~\ref{assume:sigma_alpha}, $H_c(x)\succ0$.
    \end{enumerate}
\end{lemma}
\begin{proof}
(a) Since $\Sigma_c(x)$ is symmetric, we only need to show that $\widetilde{\Sigma}_c(x)$ is symmetric.
        Using the result $\nabla \alpha_c(x)$ in Lemma~\ref{lem:misc}(c), we obtain
        \begin{equation}\label{eq:nab_eta_vcx}
            \nabla \alpha_c(x)\nabla v_c(x)^T = -\frac{1}{\delta_c\|v_c(x)\|_{W_c}}\nabla v_c(x) \alpha_c(x)\alpha_c(x)^T\nabla v_c(x)^T + \frac{\delta_c}{\|v_c(x)\|_{W_c}} \nabla v_c(x)W_c \nabla v_c(x),
        \end{equation}
        which is a symmetric matrix in $\R^{n\times n}$.

(b) We prove the result $\widetilde{\Sigma}_c(x)x=0$. 
        It holds that $\nabla v_c(x)^Tx = 2v_c(x)$ by Lemma~\ref{lem:misc}(a).
        Thus, 
        \begin{align}
            \widetilde{\Sigma}_c(x)x &= \frac{1}{2} \nabla \alpha_c(x)\nabla v_c(x)^Tx = \nabla \alpha_c(x)v_c(x) \nonumber \\
            &= -\frac{1}{\delta_c\|v_c(x)\|_{W_c}}\nabla v_c(x) \alpha_c(x)\alpha_c(x)^Tv_c(x) + \frac{\delta_c}{\|v_c(x)\|_{W_c}} \nabla v_c(x)W_c v_c(x) \nonumber \\
            &= -\nabla v_c(x)\alpha_c(x) + \nabla v_c(x)\alpha_c(x) \label{eq:null2} \\
            &= 0 \nonumber
        \end{align}
        by applying $\alpha_c(x)^Tv_c(x) = \delta_c\|v_c(x)\|_{W_c}$ and the definition of $\alpha_c(x)$~\eqref{eq:alph_x} for the equality~\eqref{eq:null2}.

(c) We show that $\widetilde{\Sigma}_c(x)$ is positive semi-definite. 
Then, $H_c(x)$ must be positive definite since it is defined as the sum 
of the positive definite matrix $\Sigma_c(x)$ and 
the positive semi-definite matrix $\widetilde{\Sigma}_c(x)$.
In fact, by definition, we have
    \begin{equation}\label{eq:tilde_sig_exp}
        \widetilde{\Sigma}_c(x) = -\frac{1}{2}\cdot\frac{\delta_c}{\|v_c(x)\|_{W_c}}\bigg(\frac{1}{\delta_c^2}\nabla v_c(x)\alpha_c(x)\alpha_c(x)^T\nabla v_c(x)^T - \nabla v_c(x)W_c\nabla v_c(x)^T\bigg).
    \end{equation}
    We consider the definiteness on the matrix inside the parenthesis of~\eqref{eq:tilde_sig_exp}. That is, for a nonzero $y\in\R^n$ we check the sign of the quantity
    \begin{equation}\label{eq:def_def}
        y^T\bigg(\frac{1}{\delta_c^2}\nabla v_c(x)\alpha_c(x)\alpha_c(x)^T\nabla v_c(x)^T - \nabla v_c(x)W_c\nabla v_c(x)^T\bigg)y.
    \end{equation}
    By the definition of $\alpha_c(x)$~\eqref{eq:alph_x} and the equation $\nabla v_c(x)$ in Lemma~\ref{lem:misc}(a), we can show that the quantity~\eqref{eq:def_def} is equal to
    \begin{align}\label{eq:def_def2}
        & \frac{4}{\|v_c(x)\|_{W_c}^2}\bigg(\sum_{i=1}^mw_c^{(i)}(x^TV_c^{(i)}x)(y^TV_c^{(i)}x)\bigg)^2 - 4\sum_{i=1}^mw_c^{(i)}(y^TV_c^{(i)}x)^2 \nonumber \\
        &= 4 \Bigg[\frac{\big(\sum_{i=1}^mw_c^{(i)}(x^TV_c^{(i)}x)(y^TV_c^{(i)}x)\big)^2}{\sum_{i=1}^mw_c^{(i)}(x^TV_c^{(i)}x)^2} - \sum_{i=1}^mw_c^{(i)}(y^TV_c^{(i)}x)^2 \Bigg].
    \end{align}
    With vectors $a_c,b_c\in\R^m$ defined as
    \begin{equation}\label{acbc}
        a_c := \begin{bmatrix}
            \sqrt{w_c^{(1)}}(x^TV_c^{(1)}x) \\
            \sqrt{w_c^{(2)}}(x^TV_c^{(2)}x) \\
            \vdots \\
            \sqrt{w_c^{(m)}}(x^TV_c^{(m)}x) \\
        \end{bmatrix},
        \qquad
        b_c := \begin{bmatrix}
            \sqrt{w_c^{(1)}}(y^TV_c^{(1)}x) \\
            \sqrt{w_c^{(2)}}(y^TV_c^{(2)}x) \\
            \vdots \\
            \sqrt{w_c^{(m)}}(y^TV_c^{(m)}x) \\
        \end{bmatrix},
    \end{equation}
    the quantity~\eqref{eq:def_def2} is further simplified as:
    \begin{equation}\label{eq:def_def3}
        4\bigg[\frac{(a_c^Tb_c)^2}{\|a_c\|_2^2} - \|b_c\|_2^2 \bigg].
    \end{equation}
    By the Cauchy-Schwarz inequality we have $(a_c^Tb_c)^2 \leq \|a_c\|_2^2\cdot\|b_c\|_2^2$ so that
    \begin{equation}\label{eq:def_def4}
        4\bigg[\frac{(a_c^Tb_c)^2}{\|a_c\|_2^2} - \|b_c\|_2^2 \bigg] \leq 4\bigg[\frac{\|a_c\|_2^2\cdot\|b_c\|_2^2}{\|a_c\|_2^2} - \|b_c\|_2^2 \bigg] = 0.
    \end{equation}
    Therefore, according to the equation of $\widetilde{\Sigma}_c(x)$ in~\eqref{eq:tilde_sig_exp}, for any nonzero $y\in\R^m$ we have 
    \[y^T\widetilde{\Sigma}_c(x)y\geq0\]
    indicating that $\widetilde{\Sigma}_c(x)$ is positive semi-definite.
\end{proof}

\begin{remark}
    For condition $\bar{c}$, Lemma~\ref{lem:Gcx}(a) and Lemma~\ref{lem:Gcx}(b) hold. However, Lemma~\ref{lem:Gcx}(c), i.e., the positive definiteness of the matrix $H_{\bar{c}}(x)$, may not hold. This is due to the negative sign present in $\alpha_{\bar{c}}(x)$, which causes the matrix $\widetilde{\Sigma}_{\bar{c}}(x)$ to be expressed as
    \[\widetilde{\Sigma}_{\bar{c}}(x) = \frac{1}{2}\cdot\frac{\delta_{\bar{c}}}{\|v_{\bar{c}}(x)\|_{W_{\bar{c}}}}\bigg(\frac{1}{\delta_{\bar{c}}^2}\nabla v_{\bar{c}}(x)\alpha_{\bar{c}}(x)\alpha_{\bar{c}}(x)^T\nabla v_{\bar{c}}(x)^T - \nabla v_{\bar{c}}(x)W_{\bar{c}}\nabla v_{\bar{c}}(x)^T\bigg),\]
    an identical expression to $\widetilde{\Sigma}_c(x)$~\eqref{eq:tilde_sig_exp} but with a sign difference.
    Consequently, $\widetilde{\Sigma}_{\bar{c}}(x)$ is negative semi-definite, which implies that the positive definiteness of $H_{\bar{c}}(x)$ is not guaranteed.
\end{remark}

\subsection{NEPv based on first order optimality conditions}

We proceed to state the first order necessary conditions of the OptNRQ~\eqref{eq:cspnrq} and derive a NEPv in Theorem~\ref{thm:first}. 


\begin{theorem}\label{thm:first}
Under Assumption~\ref{assume:sigma_alpha}, if 
$x_c$ is a local minimizer of the OptNRQ~\eqref{eq:cspnrq}, then 
$x_c$ is an eigenvector of the NEPv:
\begin{equation}\label{eq:1st_NEPv}
\begin{cases}
\Sigma_c(x)x=\lambda(\Sigma_c(x)+\Sigma_{\bar{c}}(x))x \\
x^T(\Sigma_c(x)+\Sigma_{\bar{c}}(x))x=1,
\end{cases}
\end{equation}
for some eigenvalue $\lambda>0$.
\end{theorem}
\begin{proof} 
Let us consider the constrained optimization~\eqref{eq:cspnrq_new}, an equivalent problem of the OptNRQ~\eqref{eq:cspnrq}.
    The Lagrangian function of the constrained optimization~\eqref{eq:cspnrq_new} is
    \begin{equation}\label{eq:lag_cspnrq2}
        L(x,\lambda)=x^T\Sigma_c(x)x-\lambda(x^T(\Sigma_c(x)+\Sigma_{\bar{c}}(x))x-1).
    \end{equation}
    By Lemma~\ref{lem:scx_der}, 
    the gradient of \eqref{eq:lag_cspnrq2} with respect to $x$ is
    \begin{equation}\label{eq:lag_cspnrq2_der}
        \nabla_xL(x,\lambda)=2\Sigma_c(x)x-2\lambda(\Sigma_c(x)+\Sigma_{\bar{c}}(x))x.
    \end{equation}
    Then, the first-order necessary conditions \cite[Theorem 12.1, p.321]{nocedal2006numerical} of the constrained optimization~\eqref{eq:cspnrq_new}, 
    $$
    \nabla_xL(x,\lambda)=0 
    \quad \mbox{and} \quad
    x^T(\Sigma_c(x)+\Sigma_{\bar{c}}(x))x=1,
    $$ imply that a local minimizer $x$ is an eigenvector of the NEPv~\eqref{eq:1st_NEPv}.
The eigenvalue $\lambda$ is positive due to the positive definiteness of $\Sigma_c(x)$ and $\Sigma_{\bar{c}}(x)$ under Assumption~\ref{assume:sigma_alpha}. 
\end{proof}

\paragraph{Eigenvalue ordering issue.}
Note that while a local minimizer of the OptNRQ~\eqref{eq:cspnrq} is an eigenvector of the NEPv~\eqref{eq:1st_NEPv} corresponding to {\em some} eigenvalue 
$\lambda>0$, the ordering of the desired eigenvalue $\lambda$ is unknown. 
That is, we do not know in advance which of the eigenvalues of the matrix pair $(\Sigma_c(x), \Sigma_c(x)+\Sigma_{\bar{c}}(x))$ corresponds to the desired eigenvalue $\lambda$.
To address this, we introduce an alternative NEPv for which the eigenvalue ordering of $\lambda$ is known in advance.

\subsection{NEPv based on second order optimality conditions}

Utilizing Lemma~\ref{lem:Gcx}, we derive another NEPv. The eigenvalue ordering of this alternative NEPv can be established using the second order optimality conditions of the OptNRQ~\eqref{eq:cspnrq}. 


\begin{theorem}\label{thm:second}~
    Under Assumption~\ref{assume:sigma_alpha},
    \begin{enumerate}[label=(\alph*)]

    \item If $x_c$ is a local minimizer of the OptNRQ~\eqref{eq:cspnrq}, 
    then $x_c$ is an eigenvector of the following NEPv:
    \begin{equation}\label{eq:2nd_NEPv}
        \begin{cases}
            H_c(x)x=\lambda \big(H_c(x)+H_{\bar{c}}(x)\big)x \\
            x^T\big(H_c(x)+H_{\bar{c}}(x) \big)x=1,
        \end{cases}
    \end{equation}
    corresponding to the smallest positive eigenvalue $\lambda$.
The NEPv~\eqref{eq:2nd_NEPv} is  referred to as the {OptNRQ-nepv}.
    
    \item If $\lambda$ is the smallest positive eigenvalue of the OptNRQ-nepv~\eqref{eq:2nd_NEPv} and is simple, then the corresponding eigenvector $x_c$ of the OptNRQ-nepv~\eqref{eq:2nd_NEPv} is a strict local minimizer of the OptNRQ~\eqref{eq:cspnrq}.
    \end{enumerate}
\end{theorem}

\begin{proof}  
For result (a):
Let $x$ be a local minimizer of the constrained 
optimization~\eqref{eq:cspnrq_new}, an equivalent problem of 
the OptNRQ~\eqref{eq:cspnrq}. By Theorem~\ref{thm:first}, 
$x$ is an eigenvector of the NEPv~\eqref{eq:1st_NEPv} for {\em some} 
eigenvalue $\lambda>0$. By utilizing the property 
$H_c(x)x=\Sigma_c(x)x$ in Lemma~\ref{lem:Gcx}(b), 
$x$ is also an eigenvector of the OptNRQ-nepv~\eqref{eq:2nd_NEPv} 
corresponding to the same eigenvalue $\lambda>0$.

Next, we show that $\lambda$ must be the smallest positive eigenvalue 
of the OptNRQ-nepv~\eqref{eq:2nd_NEPv} using the second order necessary 
conditions of a constrained optimization~\cite[Theorem 12.5, p.332]{nocedal2006numerical}.

Both $H_c(x)$ and $H_{\bar{c}}(x)$ are symmetric with $H_c(x)\succ0$.
Note that the positive definiteness of $H_{\bar{c}}(x)$ 
and $H_c(x)+H_{\bar{c}}(x)$ are not guaranteed.
This implies that the eigenvalues of the matrix pair 
$(H_c(x),H_c(x)+H_{\bar{c}}(x))$ are real (possibly including infinity) 
but not necessarily positive. 
We denote the eigenvalues of the OptNRQ-nepv~\eqref{eq:2nd_NEPv} 
as $\lambda_1\leq\lambda_2\leq\cdots\leq\lambda_n$ and 
the corresponding eigenvectors as $v_1,v_2,\ldots,v_n$ so that
\begin{equation}\label{eq:G_eigs}
H_c(x)v_i=\lambda_i(H_c(x)+H_{\bar{c}}(x))v_i
\quad\mbox{for}\quad i=1,\ldots,n.
\end{equation}
We take the eigenvectors to be $H_c(x)$-orthogonal, i.e.,
\begin{equation}\label{eq:G-orth}
    v_i^TH_c(x)v_j=
    \begin{cases}
        1\qquad\mbox{if}\quad i = j \\
        0\qquad\mbox{otherwise}.
    \end{cases}
\end{equation}
As $(\lambda,x)$ is an eigenpair of the OptNRQ-nepv~\eqref{eq:2nd_NEPv}, we have 
\begin{equation}\label{eq:eigenpair}
    \lambda=\lambda_{j_0}, \quad x=v_{j_0}
\end{equation}
for some $1\leq j_0\leq n$. We will deduce {the order of $\lambda$} 
from the second order necessary conditions of the constrained 
optimization~\eqref{eq:cspnrq_new} for a local minimizer $x$:
\begin{equation}\label{eq:2nd_order}
z^T\nabla_{xx}L(x,\lambda)z\geq0 
\quad\forall z\in\R^n\quad\mbox{such that}\quad 
z^T(\Sigma_c(x)+\Sigma_{\bar{c}}(x))x=0, 
\end{equation}
where $L(x,\lambda)$ is the Lagrangian function defined 
in \eqref{eq:lag_cspnrq2}. Let $z\in\R^n$ be a vector satisfying 
the condition~\eqref{eq:2nd_order}, i.e., 
$z^T(\Sigma_c(x)+\Sigma_{\bar{c}}(x))x=0$. 
We can write the expansion of $z$ in terms of the eigenbasis $\{v_i\}_{i=1}^n$:
$$z=\sum_{i=1}^na_iv_i$$
where $a_i=z^TG_c(x)v_i$ by using \eqref{eq:G-orth}. 
In particular, by \eqref{eq:G_eigs} and \eqref{eq:eigenpair}, we have 
\[ 
a_{j_0}=z^TH_c(x)v_{j_0}
=\lambda_{j_0} z^T(H_c(x)+H_{\bar{c}}(x))v_{j_0}
=\lambda z^T(\Sigma_c(x)+\Sigma_{\bar{c}}(x))x=0.
\] 
Therefore, $z$ can be written as
$$z=\sum_{i=1,i\neq {j_0}}^na_iv_i.$$ 
We use Lemma~\ref{lem:scx_der_der} to compute the gradient of 
$\nabla_xL(x,\lambda)$ in~\eqref{eq:lag_cspnrq2_der}:
$$
\nabla_{xx}L(x,\lambda)=2H_c(x)-2\lambda(H_c(x)+H_{\bar{c}}(x)).$$
So, the inequality of \eqref{eq:2nd_order} becomes 
\begin{equation}\label{eq:2nd_order_ineq}
2\cdot z^T\Big(H_c(x)-\lambda(H_c(x)+H_{\bar{c}}(x))\Big)z
=2\sum_{i=1,i\neq {j_0}}^na_i^2\bigg(1-\frac{\lambda_{j_0}}{\lambda_i}\bigg)\geq0
\end{equation}
by using~\eqref{eq:G-orth} and substituting $\lambda=\lambda_{j_0}$ from~\eqref{eq:eigenpair}.
Since \eqref{eq:2nd_order_ineq} holds for arbitrary $a_i$, the inequality~\eqref{eq:2nd_order_ineq} implies that $1-\frac{\lambda_{j_0}}{\lambda_i}\geq0,\quad\forall i=1,\ldots,n$ with $i\neq {j_0}$. As $\lambda=\lambda_{j_0}>0$, this establishes $\lambda_i\geq\lambda$ for all positive eigenvalues $\lambda_i$.

{For result (b):}
If $\lambda=\lambda_{j_0}$ is the smallest positive eigenvalue of the OptNRQ-nepv~\eqref{eq:2nd_NEPv}, and $\lambda$ is simple, then for the corresponding eigenvector $x$ of the OptNRQ-nepv~\eqref{eq:2nd_NEPv} we have
$$z^T\nabla_{xx}L(x,\lambda)z>0,\quad\forall z\in\R^n\quad\mbox{such that}\quad z^T(\Sigma_c(x)+\Sigma_{\bar{c}}(x))x=0\quad\mbox{and}\quad z\neq0$$
which is precisely the second order sufficient conditions \cite[Theorem 12.6, p.333]{nocedal2006numerical} for $x$ to be a strict local minimizer of the constrained optimization~\eqref{eq:cspnrq_new}. Consequently, $x$ is a strict local minimizer of the OptNRQ~\eqref{eq:cspnrq}.
\end{proof}

In the following example, we highlight the contrast concerning the order 
of the eigenvalue $\lambda$ with respect to the NEPv~\eqref{eq:1st_NEPv} 
and the OptNRQ-nepv~\eqref{eq:2nd_NEPv}.
Specifically, we demonstrate that while the order of $\lambda$ 
for the NEPv~\eqref{eq:1st_NEPv} cannot be predetermined, it corresponds 
to the smallest positive eigenvalue of the OptNRQ-nepv~\eqref{eq:2nd_NEPv}.

\begin{example}[Eigenvalue ordering issue] \label{eg:eigordering}
{\rm 
We solve the OptNRQ~\eqref{eq:cspnrq} for conditions $c\in\{-,+\}$ using 
a synthetic dataset where the trials are generated by 
the synthetic procedure described in Section~\ref{subsec:data}.
We generate 50 trials for each condition.
The tolerance set radius is chosen as $\delta_c=6$ for both conditions.
Let $\widehat{x}_-$ and $\widehat{x}_+$ denote the solutions obtained using the algorithm OptNRQ-nepv (discussed in Section~\ref{subsec:scf}) and verified by the Riemannian trust region (RTR) algorithm.

We first consider the NEPv~\eqref{eq:1st_NEPv}. For $\widehat{x}_-$, the first two eigenvalues of the matrix pair $(\Sigma_-(\widehat x_-),\Sigma_-(\widehat x_-)+\Sigma_+(\widehat x_-))$ are 
    \[
     \lambda_1=0.3486, \quad \underline{\lambda_2=0.4286},
    \]
    where the eigenvector corresponding to the underlined eigenvalue is $\widehat x_-$.
    Meanwhile, for $\widehat{x}_+$, the first few eigenvalues of the matrix pair $(\Sigma_+(\widehat{x}_+),\Sigma_-(\widehat{x}_+)+\Sigma_+(\widehat{x}_+))$ are 
    \[
     \lambda_1=0.3035, \quad \ldots, \quad \underline{\lambda_6=0.5120},
    \]
    where the eigenvector corresponding to the underlined eigenvalue is $\widehat x_+$.
    As we can see, the eigenvalue ordering of the NEPv~\eqref{eq:1st_NEPv} for $\widehat x_-$ is different from that of $\widehat x_+$.
    Based on this observation, we conclude that we could not have known in advance the order of the desired eigenvalue of the NEPv~\eqref{eq:1st_NEPv}. 


Next, let us consider the OptNRQ-nepv~\eqref{eq:2nd_NEPv}.
For $\widehat{x}_-$, the first two eigenvalues of the matrix pair
    $(H_-(\widehat x_-),H_-(\widehat x_-)+H_+(\widehat x_-))$ are
    \[
    \underline{\lambda_1=0.4286}, \quad \lambda_2=0.5155,
    \]
    where the eigenvector corresponding to the underlined eigenvalue is $\widehat x_-$.
Meanwhile, for $\widehat{x}_+$, the first three eigenvalues of the matrix pair
$(H_+(\widehat{x}_+),H_-(\widehat{x}_+)+H_+(\widehat{x}_+))$ are
    \[
    \lambda_1=-4.7463, \quad \underline{\lambda_2=0.5120}, \quad \lambda_3=0.5532,
    \]
    where the eigenvector corresponding to the underlined eigenvalue is $\widehat x_+$.
As we can see, for both $\widehat{x}_-$ and $\widehat{x}_+$, the desired eigenvalue is the smallest positive eigenvalue of the OptNRQ-nepv~\eqref{eq:2nd_NEPv}. 
\hfill $\Box$



}\end{example}

\begin{remark}
We do not pursue solving the NEPv~\eqref{eq:1st_NEPv} due to the challenge posed by its unknown eigenvalue ordering. 
    We will focus on solving the OptNRQ-nepv~\eqref{eq:2nd_NEPv}, for which we can always aim to compute the eigenvector associated with the smallest positive eigenvalue.
\end{remark}


\section{Algorithms}\label{sec:alg}

In this section, we consider two algorithms for solving 
the OptNRQ~\eqref{eq:cspnrq}.
The first algorithm is a fixed-point type iteration scheme applied 
directly on the OptNRQ~\eqref{eq:cspnrq}.
The second algorithm employs a self-consistent field (SCF) iteration on 
the OptNRQ-nepv~\eqref{eq:2nd_NEPv} formulation of the OptNRQ~\eqref{eq:cspnrq}.

\subsection{Fixed-point type iteration}\label{subsec:fp}
As presented in \cite{kawanabe2014robust}, 
one natural idea for solving the OptNRQ~\eqref{eq:cspnrq} is using the following 
fixed-point type iteration:
\begin{equation}\label{eq:fixed_itr}
x_c^{(k+1)}\longleftarrow\argmin_{x\neq0}
\frac{x^T\Sigma_c(x_c^{(k)})x}{x^T(\Sigma_c(x_c^{(k)})+\Sigma_{\bar{c}}(x_c^{(k)}))x}.
\end{equation}
We refer to this scheme as the {OptNRQ-fp}. In the OptNRQ-fp~\eqref{eq:fixed_itr}, a generalized Rayleigh quotient optimization involving constant matrices $\Sigma_c(x_c^{(k)})$ and $\Sigma_c(x_c^{(k)})+\Sigma_{\bar{c}}(x_c^{(k)})$ is solved to obtain the next vector $x_c^{(k+1)}$.
By the Courant-Fischer variational principle \cite{parlett1998symmetric}, 
$x_c^{(k+1)}$ of the OptNRQ-fp~\eqref{eq:fixed_itr} satisfies
the generalized eigenvalue problem 
\begin{equation}\label{eq:csp_gep2}
\Sigma_c(x_c^{(k)})\,x =\lambda \big(\Sigma_c(x_c^{(k)})+\Sigma_{\bar{c}}(x_c^{(k)})\big)\,x,
\end{equation}
where  $\lambda$ is the smallest eigenvalue. We determine its convergence by the relative residual norm of \eqref{eq:csp_gep2}.
Thus, when the OptNRQ-fp~\eqref{eq:fixed_itr} converges, we expect the solution to be an eigenvector with respect to the smallest eigenvalue of the NEPv~\eqref{eq:1st_NEPv}. 
However, as illustrated in Example~\ref{eg:eigordering}, the solution of the OptNRQ~\eqref{eq:cspnrq} may not be an eigenvector with respect to the smallest eigenvalue of the NEPv~\eqref{eq:1st_NEPv}.
Consequently, even when the OptNRQ-fp~\eqref{eq:fixed_itr} converges, its solution may not be the solution of the OptNRQ~\eqref{eq:cspnrq}.
In addition, as illustrated in Section~\ref{subsec:conv}, 
the OptNRQ-fp~\eqref{eq:fixed_itr} often fails to converge.

\subsection{Self-consistent field iteration}\label{subsec:scf}

\paragraph{Plain SCF.}
The SCF iteration, originating from solving the Kohn-Sham equation in physics and quantum chemistry \cite{cances2001self}, is a popular algorithm for solving NEPvs \cite{bai2018robust, bai2022sharp, cai2018eigenvector}. 
The basic idea of the SCF iteration is to iteratively fix the eigenvector dependence in the matrices to solve the standard eigenvalue problem.

Applied to the OptNRQ-nepv~\eqref{eq:2nd_NEPv}, the SCF iteration is given by
\begin{equation}\label{eq:SCF_map2}
x_c^{(k+1)}\longleftarrow
\mbox{an eigenvector w.r.t. } 
\lambda^+_{\min}(H_c(x_c^{(k)}),H_c(x_c^{(k)})+H_{\bar{c}}(x_c^{(k)})),
\end{equation}
where $\lambda^{+}_{\min}(A,B)$ refers to the smallest positive eigenvalue of $(A,B)$.

The converged solution of the SCF iteration~\eqref{eq:SCF_map2} is 
an eigenvector with respect to the smallest positive eigenvalue 
of the OptNRQ-nepv~\eqref{eq:2nd_NEPv}. 
Furthermore, by Theorem~\ref{thm:second}, if the corresponding eigenvalue 
is simple, then the solution is a strict local minimizer 
of the OptNRQ~\eqref{eq:cspnrq}.
Therefore, in contrast to the OptNRQ-fq~\eqref{eq:fixed_itr}, the SCF iteration~\eqref{eq:SCF_map2} provides greater confidence in its converged solution as a local minimizer of the OptNRQ~\eqref{eq:cspnrq}.

\paragraph{Line search.}
A potential issue with the plain SCF iteration~\eqref{eq:SCF_map2} is that 
the value of the objective function $q_c(x)$ of the OptNRQ~\eqref{eq:cspnrq} may not decrease monotonically. 
That is, denoting $\widetilde{x}_c^{(k+1)}$ as the next iterate of $x_c^{(k)}$ in the plain SCF iteration~\eqref{eq:SCF_map2}, it may occur that
$q_c(\widetilde{x}_c^{(k+1)})\nless q_c(x_c^{(k)})$.

To address this issue, we employ the line search to ensure the sufficient decrease of the objective function $q_c(x)$ \cite[p.21,22]{nocedal2006numerical}. 
In the line search, a search direction $d^{(k)}\in\R^n$ and a stepsize $\beta^{(k)}>0$ are determined such that the updated iterate from $x_c^{(k)}$ is given by
\begin{equation} \label{eq:xupdate}
x_c^{(k+1)} = x_c^{(k)} + \beta^{(k)} d^{(k)}
\end{equation}
to satisfy
\begin{equation}\label{eq:line}
q_c(x_c^{(k+1)}) < q_c(x_c^{(k)}).
\end{equation}
The effectiveness of the line search method relies on making appropriate choices of the search direction $d^{(k)}$ and the stepsize $\beta^{(k)}$. A desirable search direction is a descent direction of the objective function $q_c(x)$ at $x_c^{(k)}$, namely
\begin{equation}\label{eq:descent}
(d^{(k)})^T\nabla q_c(x_c^{(k)}) < 0,
\end{equation}
where, by applying Lemma~\ref{lem:scx_der}, the gradient of $q_c(x)$ is
\begin{equation}\label{eq:q_der}
    \nabla q_c(x) = \frac{2}{x^T(\Sigma_c(x)+\Sigma_{\bar{c}}(x))x}\bigg(\Sigma_c(x)-q_c(x)(\Sigma_c(x)+\Sigma_{\bar{c}}(x))\bigg) x.
\end{equation}
In the following, we discuss a strategy for choosing a descent direction $d^{(k)}$. 
In this strategy, we first consider the search direction to be
\begin{equation}\label{eq:line_dk}
d^{(k)} = \widetilde x_c^{(k+1)}-x_c^{(k)},
\end{equation}
where $\widetilde{x}_c^{(k+1)}$ is from the plain SCF iteration~\eqref{eq:SCF_map2}.
In fact, $d^{(k)}$ as defined in~\eqref{eq:line_dk} is often a descent direction except for a rare case.
To see this, note that for the direction~\eqref{eq:line_dk} we have
\begin{align*}
    \frac{1}{2}(d^{(k)})^T\nabla q_c(x_c^{(k)}) 
    &= \big(\lambda_{k+1}-q_c(x_c^{(k)})\big)\frac{(\widetilde{x}_c^{(k+1)})^T(\Sigma_c(x_c^{(k)}) + \Sigma_{\bar{c}}(x_c^{(k)}))x_c^{(k)}}{(x_c^{(k)})^T(\Sigma_c(x_c^{(k)}) + \Sigma_{\bar{c}}(x_c^{(k)}))x_c^{(k)}} \\
    &= \frac{t^{(k)}}{(x_c^{(k)})^T(\Sigma_c(x_c^{(k)}) + \Sigma_{\bar{c}}(x_c^{(k)}))x_c^{(k)}},
\end{align*}
where 
\begin{equation}\label{eq:line_tk}
t^{(k)} := (\lambda_{k+1} - q_c(x_c^{(k)})) \cdot (\widetilde{x}_c^{(k+1)})^T(\Sigma_c(x_c^{(k)}) + \Sigma_{\bar{c}}(x_c^{(k)}))x_c^{(k)}.
\end{equation}
$\lambda_{k+1}$ corresponds to $\lambda^+_{\min}(H_c(x_c^{(k)}), H_c(x_c^{(k)})+H_{\bar{c}}(x_c^{(k)}))$. Since the term $(x_c^{(k)})^T(\Sigma_c(x_c^{(k)}) + \Sigma_{\bar{c}}(x_c^{(k)}))x_c^{(k)}$ is always positive, the sign of $t^{(k)}$ determines whether $d^{(k)}$~\eqref{eq:line_dk} is a descent direction.
\begin{itemize}

    \item $t^{(k)}$ is negative. In this case, $d^{(k)}=\widetilde x_c^{(k+1)}-x_c^{(k)}$ is a descent direction.

    \item $t^{(k)}$ is positive. In this case, we can simply scale $\widetilde x_c^{(k+1)}$ by $-1$, i.e., 
    \begin{equation}\label{eq:flip}
        \widetilde x_c^{(k+1)}\longleftarrow -\widetilde x_c^{(k+1)},
    \end{equation}
    to ensure that $d^{(k)}=\widetilde x_c^{(k+1)}-x_c^{(k)}$ is a descent direction. 

    \item $t^{(k)}$ is zero. We can assume that $\lambda_{k+1}-q_c(x_c^{(k)})\neq0$, as otherwise, it would indicate that the iteration has already converged. Thus, only when $(\widetilde{x}_c^{(k+1)})^T(\Sigma_c(x_c^{(k)}) + \Sigma_{\bar{c}}(x_c^{(k)}))x_c^{(k)}=0$, we have that $t^{(k)}=0$. In this rare case, we set the search direction as the steepest-descent direction of the objective function:
\begin{equation}\label{eq:steepest_d}
        d^{(k)} = -\frac{\nabla q_c(x_c^{(k)})}{\|\nabla q_c(x_c^{(k)})\|^2}.
    \end{equation}

\end{itemize}
By adopting this strategy, we ensure that we can always choose a descent direction. 
Moreover, we scarcely choose the steepest-descent direction~\eqref{eq:steepest_d} and avoid replicating the steepest descent method, which can suffer from slow convergence \cite{nocedal2006numerical}.

For determining the stepsize $\beta^{(k)}$, we follow the standard line search practice where the stepsize is computed to satisfy the following sufficient decrease condition (also known as the Armijo condition):
\begin{equation}\label{eq:armijo}
q_c(x_c^{(k)}+\beta^{(k)} d^{(k)}) \leq q_c(x_c^{(k)}) + \mu\beta^{(k)} \cdot (d^{(k)})^T\nabla q_c(x_c^{(k)}),
\end{equation}
for some $\mu\in(0,1)$. 
As the objective function $q_c(x)$ is uniformly bounded below by zero, there exists a $\beta^{(k)}$ satisfying the sufficient decrease condition~\eqref{eq:armijo} \cite[Lemma 3.1]{nocedal2006numerical}.

To obtain such a stepsize $\beta^{(k)}$, we apply the Armijo backtracking \cite[Algorithm 3.1, p.37]{nocedal2006numerical}, which ensures that the stepsize $\beta^{(k)}$ is short enough to have a sufficient decrease but not too short. 
We start with an initial $\beta^{(k)}=1$ and update the step length as $\beta^{(k)}\longleftarrow\tau\beta^{(k)}$, for some $\tau\in(0,1)$, until the sufficient decrease condition~\eqref{eq:armijo} is satisfied.

\paragraph{Algorithm pseudocode.}
The algorithm for the SCF iteration~\eqref{eq:SCF_map2} with line search is shown in Algorithm~\ref{alg:SCF}. We refer to Algorithm~\ref{alg:SCF} as the {OptNRQ-nepv}. A few remarks are in order. 

\begin{algorithm}[t]
\caption{OptNRQ-nepv}\label{alg:SCF}
\begin{flushleft}
\textbf{Input}: Initial $x_c^{(0)}\in\R^n$, line search factors $\mu,\tau\in(0,1)$ (e.g. $\mu=\tau=0.01$), tolerance $tol$. \\
\textbf{Output}: Approximate solution $\widehat{x}_c$ to the OptNRQ~\eqref{eq:cspnrq}.
\end{flushleft}
\begin{algorithmic}[1]
    \FOR{$k=0,1,\ldots$}
        \STATE Set $A_k=H_c(x_c^{(k)})$, $B_k=H_c(x_c^{(k)})+H_{\bar{c}}(x_c^{(k)})$, and $q_k=\frac{(x_c^{(k)})^TA_kx_c^{(k)}}{(x_c^{(k)})^TB_kx_c^{(k)}}$.
        \IF{$\frac{||A_kx_c^{(k)}-q_kB_kx_c^{(k)}||_2}{||A_kx_c^{(k)}||_2+q_k||B_kx_c^{(k)}||_2}<\mbox{tol}$}
            \STATE return $\widehat{x}_c=x_c^{(k)}$.
        \ENDIF
        \STATE Compute the smallest positive eigenpair 
$(\lambda^{(k+1)},x_c^{(k+1)})$ of $(A_k,B_k)$.
        \IF{$q_c(x_c^{(k+1)})\geq q_c(x_c^{(k)})$}
            \STATE Set $t^{(k)} = (\lambda^{(k+1)} - q_c(x_c^{(k)})) \cdot (x_c^{(k+1)})^TB_kx_c^{(k)}$.
            \IF{$|t^{(k)}|<\mbox{tol}$}
                \STATE Set $d^{(k)}=-\frac{\nabla q_c(x_c^{(k)})}{\|\nabla q_c(x_c^{(k)})\|^2}$.
            \ELSIF{$t^{(k)}>0$}
                \STATE Set $d^{(k)}=-x_c^{(k+1)}-x_c^{(k)}$.
            \ELSE
                \STATE Set $d^{(k)}=x_c^{(k+1)}-x_c^{(k)}$.
            \ENDIF
            \STATE Set $\beta^{(k)}=1$.
            \WHILE{$q_c(x_c^{(k)}+\beta^{(k)} d^{(k)}) > q_c(x_c^{(k)}) + \mu\beta^{(k)} \cdot (d^{(k)})^T\nabla q_c(x_c^{(k)})$}
                \STATE{$\beta^{(k)}:=\tau\beta^{(k)}$}.
            \ENDWHILE
            \STATE $x_c^{(k+1)}:=x_c^{(k)}+\beta^{(k)} d^{(k)}$.
        \ENDIF
    \ENDFOR
\end{algorithmic}
\end{algorithm}



\begin{itemize} 


\item In line 3, the stopping criterion checks whether the OptNRQ-nepv~\eqref{eq:2nd_NEPv} is satisfied.

\item The smallest positive eigenpair in line 6 can be efficiently computed using sparse eigensolvers such as the implicitly restarted Arnoldi method \cite{sorensen1997implicitly} or Krylov-Schur algorithm \cite{stewart2002krylov}.

\item If the value of the objective function $q_c(x)$ does not decrease at the next iterate (line 7), the algorithm takes appropriate steps (lines 8 to 15) to find a descent direction, and the Armijo backtracking (lines 16 to 19) is used to obtain an appropriate step length $\beta^{(k)}$, and consequently, the next iterate is updated in line 20.


\end{itemize}

\begin{remark}
While a formal convergence proof of the OptNRQ-nepv (Algorithm~\ref{alg:SCF}) is still subject to further study, our numerical experiments demonstrate both convergence and quadratic local convergence. 
For the proof of quadratic local convergence for the SCF iteration pertaining to a general robust Rayleigh quotient optimization, we refer the reader to \cite{bai2018robust}.
\end{remark}


\section{Numerical Experiments}\label{sec:num}

We first provide descriptions of the synthetic datasets and 
real-world EEG datasets to be used in the experiments. 
The specific experimental settings are outlined in the second subsection.
In the third subsection, we present the convergence behavior and running time of the OptNRQ-nepv (Algorithm~\ref{alg:SCF}), the OptNRQ-fp~\eqref{eq:fixed_itr}, and manifold optimization algorithms. 
Lastly, we present the classification results obtained using the standard CSP, 
the OptNRQ-fp and the OptNRQ-nepv for both the synthetic and real-world datasets.

\subsection{Datasets}\label{subsec:data}

\paragraph{Synthetic dataset.}
The synthetic dataset, utilized in \cite{kawanabe2014robust,bai2018robust}, comprises of trials 
\[ 
Y_c^{(i)}=\begin{bmatrix} x_1 & x_2 & \cdots x_t \end{bmatrix} \in\R^{n\times t},
\]
where, for $i=1,2,\ldots,t$, $x_i\in\R^n$ is a synthetic signal of $n$ channels. The synthetic signals are randomly generated using the following linear mixing model:
\begin{equation}\label{eq:lin_mix_model}
    x_i = A\begin{bmatrix}s^{n_1}\\s^{n_2}\end{bmatrix} + \epsilon.
\end{equation}
The matrix $A\in\R^{n\times n}$ is a random rotation matrix that simulates the smearing effect that occurs when the source signals are transmitted through the skull and the skin of the head and recorded by multiple electrodes on the scalp during EEG recordings. 
The discriminative sources $s^{n_1}\in\R^{n_1}$ represent the brain sources that exhibit changes in rhythmic activity during motor imagery, whereas the nondiscriminative sources $s^{n_2}\in\R^{n_2}$ represent the brain sources that maintain their rhythmic activity throughout. 
The nonstationary noise $\epsilon\in\R^n$ corresponds to artifacts that occur during the experiment, such as blinking and muscle movements. 
We consistently employ the same procedure to generate trials 
when using the synthetic dataset.
We refer to this procedure simply as the synthetic procedure.
In the synthetic procedure, the signals in a trial consist of $n=10$ channels. 
For either condition, 2 sources of the signals are discriminative while the other 8 sources are nondiscriminative.
The first discriminative source is sampled from $\mathcal{N}(0, 0.2)$\footnote{$\mathcal{N}(\mu, \sigma^2)$ denotes the Gaussian distribution where $\mu$ is the mean and $\sigma^2$ is the variance.} for condition `$-$' and $\mathcal{N}(0, 1.8)$ for condition `$+$'.
The second discriminative source is sampled from $\mathcal{N}(0, 1.4)$ for condition `$-$' and $\mathcal{N}(0, 0.6)$ for condition `$+$'.
For either condition, all 8 nondiscriminative sources are sampled from $\mathcal{N}(0,1)$ and the 10 nonstationary noise sources of $\epsilon$ are sampled from $\mathcal{N}(0,2)$.
Finally, there are $t=200$ many signals in a trial.

\paragraph{Dataset IIIa.}
The dataset is from the BCI competition III \cite{blankertz2006bci} and includes EEG recordings from three subjects performing left hand, right hand, foot, and tongue motor imageries. 
EEG signals were recorded from 60 electrode channels, and each subject had either six or nine runs. 
Within each run, there were ten trials of motor imagery performed for each motor condition. 
For more detailed information about the dataset, please refer to the competition's website.\footnote{\href{http://www.bbci.de/competition/iii/desc_IIIa.pdf}{http://www.bbci.de/competition/iii/desc\_IIIa.pdf}}

We use the recordings corresponding to left hand and right hand motor imageries, and represent them as conditions `$-$' and `$+$', respectively.
Following the preprocessing steps suggested in \cite{samek2012stationary,lotte2010regularizing}, we apply the 5th order Butterworth filter to bandpass filter the EEG signals in the 8-30Hz frequency band. 
We then extract the features of the EEG signals from the time segment from 0s to 3.0s after the cue instructing the subject to perform motor imagery. 
In summary, for each subject, there are $N_c=60\mbox{ or }90$ trials $Y_c^{(i)}\in\R^{60\times750}$ for both conditions $c\in\{-,+\}$.

\paragraph{Distraction dataset.}
The distraction dataset is obtained from the BCI experiment \cite{brandl2015bringing,brandl2016brain,brandl2020motor} that aimed to simulate real-world scenarios by introducing six types of distractions on top of the primary motor imagery tasks (left hand and right hand). 
The six secondary tasks included no distraction (Clean), closing eyes (Eyes), listening to news (News), searching the room for a particular number (Numbers), watching a flickering video (Flicker), and dealing with vibro-tactile stimulation (Stimulation). 
The experiment involved 16 subjects and EEG recordings were collected from 63 channels. 
The subjects performed 7 runs of EEG recordings, with each run consisting of 36 trials of left hand motor imagery and 36 trials of right hand motor imagery. 
In the first run, motor imageries were performed without any distractions, and in the subsequent six runs, each distraction was added, one at a time. 
The dataset is publicly available.\footnote{\href{https://depositonce.tu-berlin.de/handle/11303/10934.2}{https://depositonce.tu-berlin.de/handle/11303/10934.2}}

We let condition `$-$' represent the left hand motor imagery and `$+$' represent the right hand motor imagery.
For preprocessing the EEG signals, we use the provided subject-specific frequency bands and time intervals to bandpass filter and select the time segment of the EEG signals. 
In summary, there are $N_c=252$ trials $Y_c^{(i)}\in\R^{63\times t}$ for both conditions $c\in\{-,+\}$, where the number of time samples $t$ differs for each subject and ranges between $1,550$ to $3,532$.

\paragraph{Berlin dataset.}
This dataset is from the Berlin BCI experiment described in \cite{sannelli2019large,blankertz2010neurophysiological}, which involved 80 subjects performing left hand, right hand, and right foot motor imagery while EEG signals were recorded from 119 channels. 
The dataset consists of three runs of motor imagery calibration and three runs of motor imagery feedback, with each calibration run comprising 75 trials per condition and each feedback run comprising 150 trials per condition. 
Although the dataset is publicly available,\footnote{\href{https://depositonce.tu-berlin.de/handle/11303/8979?mode=simple}{https://depositonce.tu-berlin.de/handle/11303/8979?mode=simple}} only the signals from 40 subjects are provided.

We focus on the recordings corresponding to the left hand and right hand, and represent them as conditions `$-$' and `$+$', respectively.
We preprocess the EEG signals by bandpass filtering them with the 2nd order Butterworth filter in the 8-30Hz frequency band and selecting 86 out of the 119 electrode channels that densely cover the motor cortex. 
For each subject, there are $N_c=675$ trials $Y_c^{(i)}\in\R^{86\times301}$ for both conditions $c\in\{-,+\}$.

\paragraph{Gwangju dataset.}
The EEG signals in this dataset\footnote{\href{http://gigadb.org/dataset/view/id/100295/}{http://gigadb.org/dataset/view/id/100295/}} were recorded during the motor-imagery based Gwangju BCI experiment \cite{cho2017eeg}. 
EEG signals were collected from 64 channels while 52 subjects performed left and right hand motor imagery. 
Each subject participated in either 5 or 6 runs, with each run consisting of 20 trials per motor imagery condition. 
Some of the trials were labeled as `bad trials' either if their voltage magnitude exceeded a particular threshold or if they exhibited a high correlation with electromyography (EMG) activity.

The left hand and right hand motor imageries are represented as conditions `$-$' and `$+$', respectively.
Following the approach in \cite{cho2017eeg}, we discard the bad trials for each subject and exclude the subjects `s29' and `s34', who had more than 90\% of their trials declared as bad trials. 
To preprocess the data, we apply the 4th order Butterworth bandpass filter to the signals in the 8-30Hz frequency range, and select the time interval from 0.5s to 2.5s after the cue instructing the subject to perform motor imagery. 
For each subject, there are $N_c$ trials $Y_c^{(i)}\in\R^{64\times1025}$ for both conditions $c\in\{-,+\}$, where the number of trials $N_c$ differs for each subject and ranges between 83 to 240.

\paragraph{Summary.} 
In Table~\ref{tab:summary}, we summarize the datasets, which include the number of channels and time samples of each trial, as well as the number of trials and participating subjects. Among the available motor imagery types, we italicize the ones we utilize.

\begin{table}[H]
\centering
\caption{Summary of datasets.} \label{tab:summary}
\resizebox{\textwidth}{!}{
\begin{tabular}{c|ccccc} \hline
Name & {Channels}($n$) & {Time Samples} ($t$) & {Trials} ($N_c$) & {Subjects} & {Motor Imagery Type} \\ \hline
\textbf{Synthetic} & $n$ & $t$ & $N_c$ & N/A & N/A \\ \hline
\textbf{Dataset IIIa} & 60 & 750 & 60 or 90 & 3 & \begin{tabular}[c]{@{}c@{}}\textit{left hand, right hand},\\ foot, tongue\end{tabular} \\ \hline
\textbf{Distraction} & 63 & range of 1550 to 3532 & 252 & 16 & \textit{left hand, right hand} \\ \hline
\textbf{Berlin} & \begin{tabular}[c]{@{}c@{}}119\\ (86 chosen)\end{tabular} & 301 & 675 & \begin{tabular}[c]{@{}c@{}}80\\ (40 available)\end{tabular} & \begin{tabular}[c]{@{}c@{}}\textit{left hand, right hand},\\ right foot\end{tabular} \\ \hline
\textbf{Gwangju} & 64 & 1025 & range of 83 to 240 & \begin{tabular}[c]{@{}c@{}}52\\ (50 used)\end{tabular} & \textit{left hand, right hand} \\ \hline
\end{tabular}
}
\end{table}


\subsection{Experiment setting}
We utilize the BBCI toolbox\footnote{\url{https://github.com/bbci/bbci_public}} 
to preprocess the raw EEG signals.
For both conditions $c\in\{-,+\}$, we choose $m=10$ as the number of interpolation matrices $\{V_c^{(i)}\}_{i=1}^m$ and weights $\{w_c^{(i)}\}_{i=1}^m$ in the formulations of the OptNRQ~\eqref{eq:cspnrq}.

Three different approaches are considered to solve the OptNRQ~\eqref{eq:cspnrq}:
\begin{itemize}
\item 
The OptNRQ-fp~\eqref{eq:fixed_itr} that solves 
the OptNRQ~\eqref{eq:cspnrq} directly.

\item The OptNRQ-nepv (Algorithm~\ref{alg:SCF}) that solves the NEPv formulation, the OptNRQ-nepv~\eqref{eq:2nd_NEPv}, of the OptNRQ~\eqref{eq:cspnrq}.

\item Manifold optimization algorithms in Manopt~\cite{manopt} to 
solve the equivalent constrained optimization~\eqref{eq:cspnrq_man} of 
the OptNRQ~\eqref{eq:cspnrq}.
\end{itemize}
For all algorithms, the initializer is set as the solution of the CSP, i.e., 
the generalized Rayleigh quotient optimizations~\eqref{eq:gRq}. 
The tolerance \textit{tol} for the stopping criteria is set to $10^{-8}$.

The experiments were conducted using MATLAB on an Apple M1 
Pro processor with 16GB of RAM.
To advocate for reproducible research, we share the preprocessed datasets as well as our MATLAB codes for algorithms and numerical experiments presented in this paper.\footnote{Github page: \url{https://github.com/gnodking7/Robust-CSP.git}}

\subsection{Convergence behavior and timing}~\label{subsec:conv}
In this section, we examine the convergence behavior of the OptNRQ-nepv (Algorithm~\ref{alg:SCF}) in comparison to those of the OptNRQ-fp~\eqref{eq:fixed_itr} and manifold optimization algorithms. 
Our findings reveal that the OptNRQ-fp fails to converge and exhibits oscillatory behavior. 
The manifold optimization algorithms do converge, but they are relatively slow when compared to the OptNRQ-nepv. 
On the other hand, the OptNRQ-nepv converges within a small number of iterations and displays a local quadratic convergence.

\begin{example}[Synthetic dataset]\label{eg:conv_synth1}
{\rm
For condition `$-$', 50 trials are generated by following the synthetic procedure.
The tolerance set radius is chosen as $\delta_-=6$.

Four algorithms are used to solve the OptNRQ~\eqref{eq:cspnrq}:
    \begin{itemize}
    
        \item The OptNRQ-fp~\eqref{eq:fixed_itr}.
        
        \item The OptNRQ-nepv (Algorithm~\ref{alg:SCF}).

        \item The Riemannian conjugate gradient (RCG) algorithm in Manopt.

        \item The Riemannian trust region (RTR) algorithm in Manopt.
        
    \end{itemize}
    
The convergence behaviors of the algorithms are depicted in Figure~\ref{fig:conv}.
Panel (a) displays the value of the objective function $q_-(x)$ of the OptNRQ~\eqref{eq:cspnrq} at each iterate $x_-^{(k)}$, while panel (b) shows the errors, measured as the difference between $q_-(x_-^{(k)})$ and $q_-(\widehat{x}_-)$, where $\widehat{x}_-$ is the solution computed by the OptNRQ-nepv.

These results demonstrate that the OptNRQ-fp~\eqref{eq:fixed_itr} oscillates 
between two points and fails to converge.  
Meanwhile, the RCG takes 497 iterations to converge in 1.416 seconds.
The RTR algorithm takes 19 iterations to converge in 0.328 seconds.
The OptNRQ-nepv achieves convergence in just 9 iterations with 
3 line searches in 0.028 seconds.

\medskip
\begin{center} 
\begin{tabular}{c|cccc}
 & {OptNRQ-fp} & {OptNRQ-nepv} & {RCG} & {RTR} \\ \hline
{Iteration} & Does not converge & 9 (3) & 497 & 19 \\
{Time (seconds)} & 0.05 & 0.028 & 1.416 & 0.328
\end{tabular} 
\end{center} 

\medskip

Although the RTR algorithm demonstrates faster convergence compared 
to the RCG algorithm due to its utilization of second-order gradient information,
it falls short in terms of convergence speed when compared 
to the OptNRQ-nepv. 
We attribute the swift convergence of the OptNRQ-nepv to its utilization of the explicit equation associated with the second order derivative information. 
Namely, the OptNRQ-nepv takes advantage of solving the OptNRQ-nepv~\eqref{eq:2nd_NEPv} derived from the second order optimality conditions of the OptNRQ~\eqref{eq:cspnrq}.

\begin{figure}[tbhp]
    \centering
    \subfloat[]
    {{\includegraphics[scale=0.18]{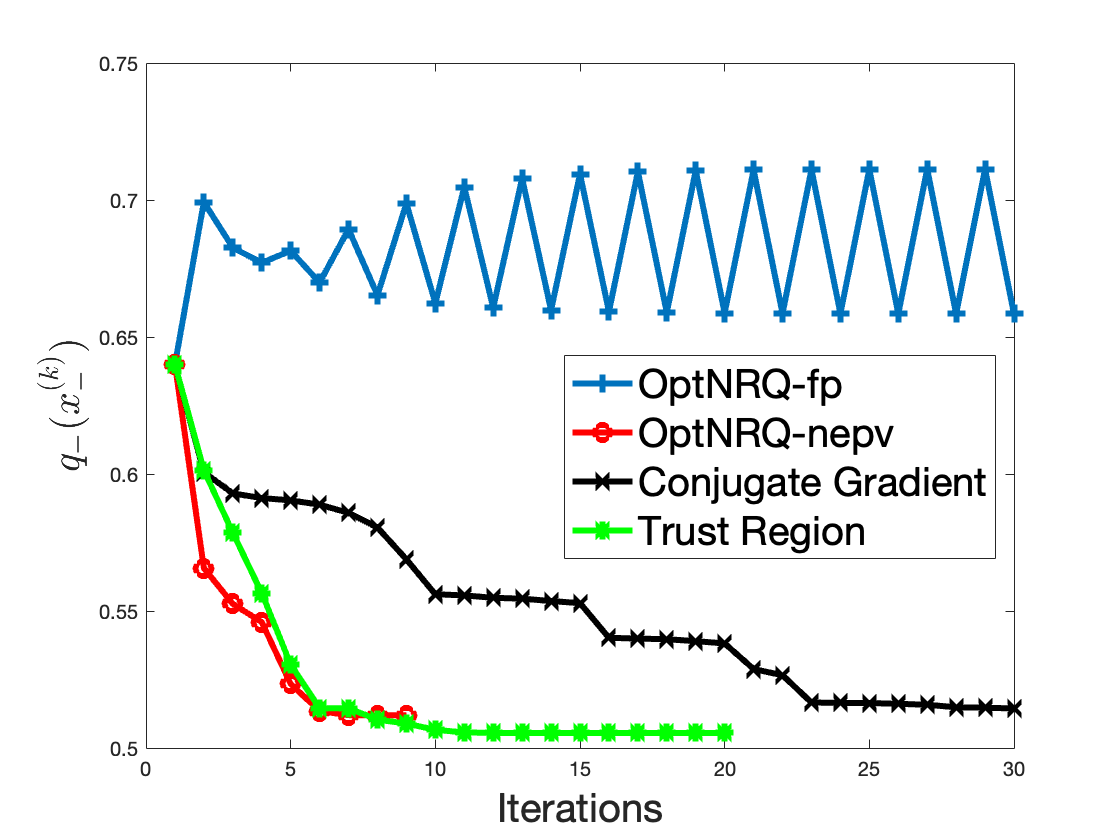} }}%
    \subfloat[]{{\includegraphics[scale=0.18]{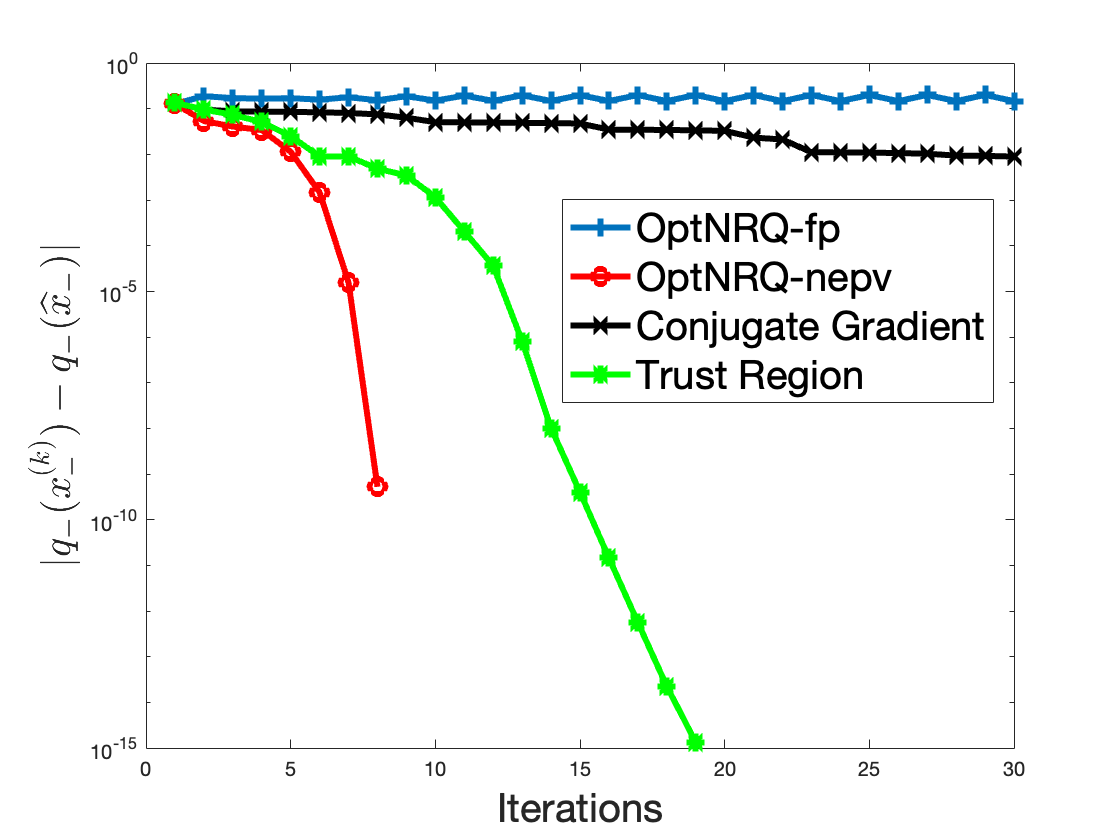} }}%
    \caption{Synthetic dataset: convergence behaviors of the OptNRQ-fp~\eqref{eq:fixed_itr}, the OptNRQ-nepv (Algorithm~\ref{alg:SCF}), the Riemannian conjugate gradient algorithm in Manopt, and the Riemannian trust region algorithm in Manopt. (a) displays the objective value $q_-(x_-^{(k)})$ of the OptNRQ~\eqref{eq:cspnrq} at each iteration, and (b) shows the errors between the objective values of the iterate $x_-^{(k)}$ and the solution $\widehat{x}_-$ obtained by the OptNRQ-nepv.}\label{fig:conv}
\end{figure}

} \end{example}

\begin{example}[Berlin dataset]\label{eg:conv_berlin1}
{\rm
We continue the convergence behavior experiments on the Berlin dataset. 
We focus on the subject \#40 and utilize the trials from the calibration runs. 
We consider the left hand motor imagery, represented by $c=-$, and set $\delta_-=0.4$ as the tolerance set radius.

Four algorithms are considered for solving the OptNRQ~\eqref{eq:cspnrq}: the OptNRQ-fp~\eqref{eq:fixed_itr}, 
the OptNRQ-nepv (Algorithm~\ref{alg:SCF}), the RCG, and the RTR.
The convergence behaviors of algorithms for this subject are shown in Figure~\ref{fig:Berlin_conv} where its two panels show analogous plots to Figure~\ref{fig:conv}.

\begin{figure}[tbhp]
    \centering
    \subfloat[]
    {{\includegraphics[scale=0.18]{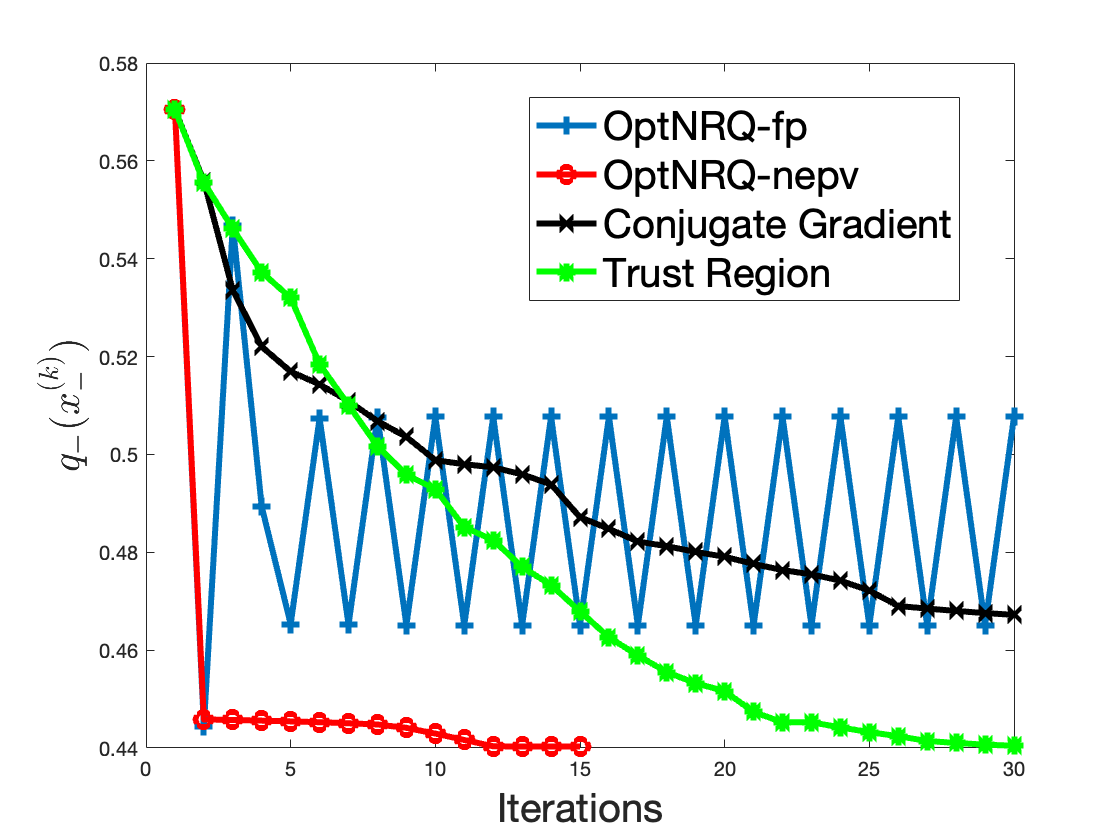} }}%
    \subfloat[]{{\includegraphics[scale=0.18]{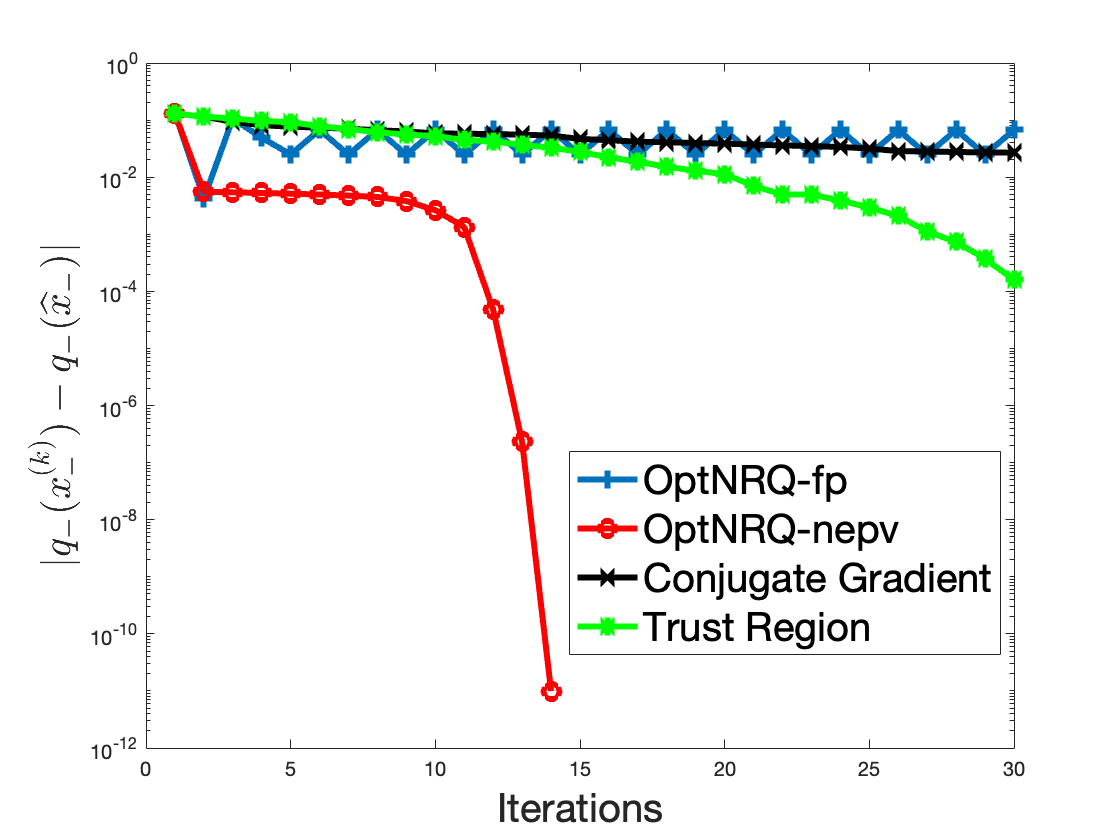} }}%
\caption{
Berlin dataset: convergence behaviors.
}\label{fig:Berlin_conv}
\end{figure}

The OptNRQ-fp~\eqref{eq:fixed_itr}, once again, exhibits oscillatory behavior and fails to converge.
The RCG algorithm fails to converge within the preset 1000 iterations and takes 38.2 seconds.
The RTR algorithm takes 73 iterations to converge in 73.792 seconds.
The OptNRQ-nepv achieves the fastest convergence, converging in 14 iterations with 4 line searches in 0.389 seconds.
Just like in Figure~\ref{fig:conv}(b), we observe a local quadratic convergence of the OptNRQ-nepv in Figure~\ref{fig:Berlin_conv}(b).

\medskip
\begin{center} 
\begin{tabular}{c|cccc}
& {OptNRQ-fp} & {OptNRQ-nepv} & {RCG} & {RTR} \\ \hline
{Iteration} & Does not converge & 14 (4) & 1000+ & 73 \\
Time (seconds) & 0.91 & 0.389 & 38.2 & 73.792
\end{tabular}   
\end{center} 
}\end{example}

\subsection{Classification results}

We follow the common signal classification procedure in BCI-EEG experiments:

\begin{enumerate}[label=(\roman*)]

    \item The principal spatial filters $x_-$ and $x_+$ are computed using the trials corresponding to the training dataset.

    \item A linear classifier is trained using the principal spatial filters and the training trials.
    
    \item The conditions of the testing trials are determined using the trained linear classifier.
        
\end{enumerate}

\paragraph{Linear classifier and classification rate.}
After computing the principal spatial filters $x_-$ and $x_+$, a common classification practice is to apply a linear classifier to the log-variance features of the filtered signals \cite{blankertz2007optimizing,blankertz2010neurophysiological,kawanabe2014robust}. 
While other classifiers that do not involve the application of the logarithm are discussed in \cite{tomioka2006logistic, tomioka2007classifying}, we adhere to the prevailing convention of using a linear classifier on the log-variance features.

For a given trial $Y$, the log-variance features $f(Y) \in \mathbb{R}^{2}$ are computed as
\begin{equation}\label{eq:log-var}
    f(Y) =
    \begin{bmatrix}
        \log(x_-^TYY^Tx_-) \\
        \log(x_+^TYY^Tx_+) \\ 
    \end{bmatrix}
\end{equation}
and a linear classifier is defined as
\begin{equation}\label{eq:classifier}
    \varphi(Y) := a^Tf(Y)-b
\end{equation}
where the sign of the classifier $\varphi(Y)$ determines the condition of $Y$. 
The weights $a\in\R^{2}$ and $b\in\R$ are determined from the training trials $\{Y_c^{(i)}\}_{i=1}^{N_c}$ using Fisher's linear discriminant analysis (LDA). Specifically, denoting $f_c^{(i)} := f(Y_c^{(i)})$ as the log-variance features, and 
\begin{equation}\label{eq:scatter}
    C_c=\sum_{i=1}^{N_c}(f_c^{(i)}-m_c)(f_c^{(i)}-m_c)^T\quad\mbox{with}\quad m_c=\frac{1}{N_c}\sum_{i=1}^{N_c}f_c^{(i)}
\end{equation}
as the scatter matrices for $c\in\{-,+\}$, the weights are determined by
\begin{equation}\label{eq:lda_weights}
a=\frac{\widetilde{a}}{\|\widetilde{a}\|_2}\quad\mbox{with}\quad 
\widetilde{a}=(C_-+C_+)^{-1}(m_--m_+),
\end{equation}
and 
\begin{equation}\label{eq:lda_weights_2}
b=\frac{1}{2}a^T(m_-+m_+).
\end{equation}

With the weights $a$ and $b$ computed as \eqref{eq:lda_weights} and \eqref{eq:lda_weights_2}, respectively, we note that the linear classifier $\varphi(Y)$~\eqref{eq:classifier} 
is positive if and only if $a^Tf(Y)$ is closer to $a^Tm_-$ than 
it is to $a^Tm_+$.
Consequently, we count a testing trial 
$Y_-^{(i)}$ in condition `$-$' to be correctly classified if $\varphi(Y_-^{(i)}) > 0$, and similarly,
a testing trial $Y_+^{(i)}$ 
in condition `$+$' to be correctly classified if $\varphi(Y_+^{(i)}) < 0$. 
Thus, using the testing trials $\{Y_-^{(i)}\}_{i=1}^{N_-}$ and $\{Y_+^{(i)}\}_{i=1}^{N_+}$, 
the {\bf classification rate} of a subject is defined as:
\begin{equation} \label{eq:classrate} 
\mbox{Classification Rate } := 
\frac{|\{\varphi(Y_-^{(i)})>0\}_{i=1}^{N_-}|}{N_-} + 
\frac{|\{\varphi(Y_+^{(i)})<0\}_{i=1}^{N_+}|}{N_+},
\end{equation}
where $|\{\varphi(Y_-^{(i)})>0\}_{i=1}^{N_-}|$ represents the number 
of times that $\varphi(Y_-^{(i)})>0$, and 
$|\{\varphi(Y_+^{(i)})<0\}_{i=1}^{N_+}|$ is defined similarly.

\begin{example}[Synthetic dataset]
{\rm 
For both conditions $c\in\{-,+\}$, we generate 50 training trials following the synthetic procedure. 
50 testing trials are also generated by the synthetic procedure for both conditions, except that the 10 nonstationary noise sources of $\epsilon$ are sampled from $\mathcal{N}(0,30)$.

We compute three sets of principal spatial filters $x_-$ and $x_+$: 
\begin{itemize}
\item[(i)] One set obtained from the CSP by solving the generalized Rayleigh quotient optimization~\eqref{eq:gRq}, 
\item[(ii)]
Another set obtained from the OptNRQ-fp~\eqref{eq:fixed_itr}, 
\item[(iii)]
The last set obtained from the OptNRQ-nepv (Algorithm~\ref{alg:SCF}). 
\end{itemize} 
We consider various tolerance set radii $\delta_c\in\{0.5,1,2,4,6,8\}$ for the OptNRQ~\eqref{eq:cspnrq}. 
When the OptNRQ-fp~\eqref{eq:fixed_itr} fails to converge, we select the solutions with the smallest objective value among its iterations.
The experiment is repeated for 100 random generations of the synthetic dataset.



The resulting classification rates are summarized as boxplots in Figure~\ref{fig:boxplots}. 
In both panels (a) and (b), the first boxplot corresponds to the classification rate of the CSP. 
The rest of the boxplots correspond to the classification rate results by the OptNRQ-fp in panel (a), and by the OptNRQ-nepv in panel (b) for different $\delta_c$ values.

\begin{figure}[tbhp]
\centering
\subfloat[]{{\includegraphics[scale=0.18]{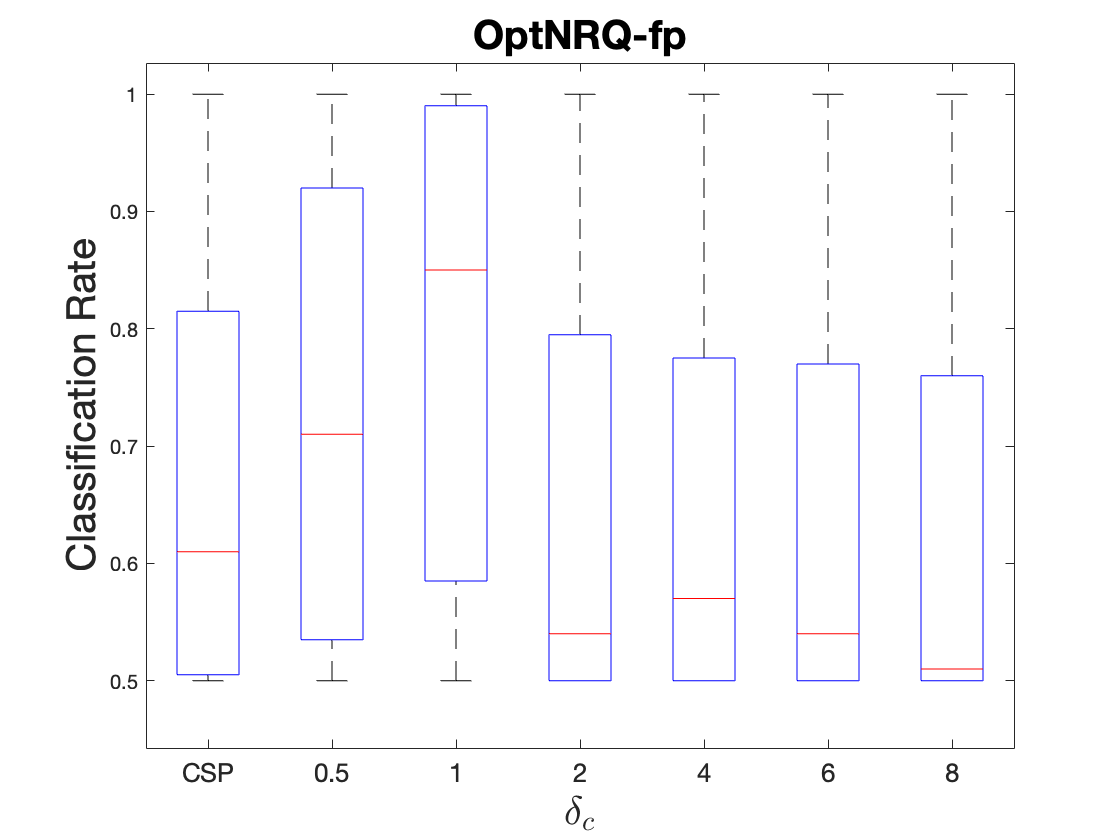} }}%
\subfloat[]{{\includegraphics[scale=0.18]{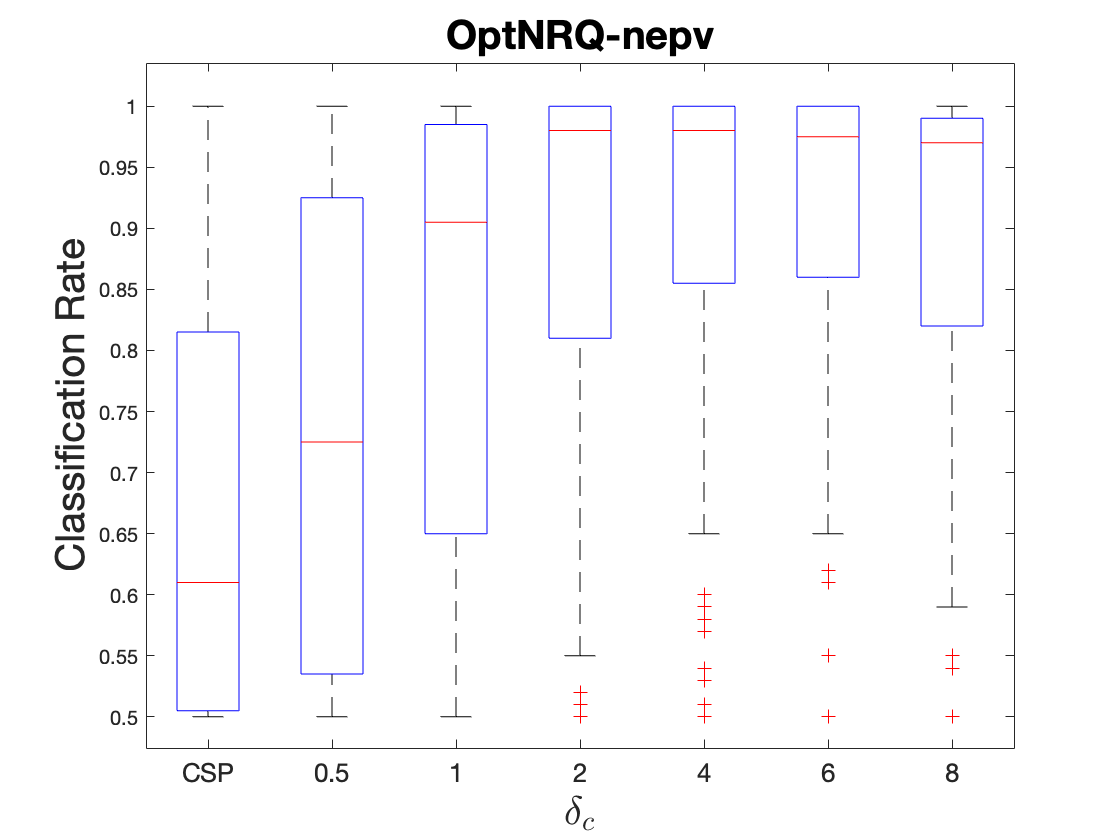} }}%
\caption{Synthetic dataset: classification rate boxplots. The first boxplot of both panel (a) and panel (b) corresponds to the classification rate of the CSP. The rest of the boxplots correspond to the classification rates of the OptNRQ-fp~\eqref{eq:fixed_itr} in panel (a), and of the OptNRQ-nepv (Algorithm~\ref{alg:SCF}) in panel (b) for the tolerance set radius $\delta_c\in\{0.5,1,2,4,6,8\}$.}\label{fig:boxplots}
\end{figure}

For small $\delta_c$ values of 0.5 and 1, the OptNRQ-fp~\eqref{eq:fixed_itr} converges, but it fails to converge for higher $\delta_c$ values. 
We observe that the performance of the OptNRQ-fp deteriorates significantly for these higher $\delta_c$ values and actually performs worse than the CSP.
In contrast, the OptNRQ-nepv demonstrates a steady increase in performance as $\delta_c$ increases beyond 1. 
Our observations show that at $\delta_c=4$ and $\delta_c=6$ the OptNRQ-nepv performs the best.
Not only does the OptNRQ-nepv improve the CSP significantly at these radius values, but it also outperforms the OptNRQ-fq by a large margin.

Additionally, we provide the iteration profile of the OptNRQ-nepv with respect to $\delta_c$ in Table~\ref{tab:itr_num}.
For each $\delta_c$, the first number represents the number of iterations, while the second number corresponds to the number of line searches.
We observe that for small values of $\delta_c$, the number of iterations is typically low, and no line searches are required.
However, as $\delta_c$ increases, both the number of iterations and line searches also increase.

\begin{table}[t]
\begin{tabular}{|c|cc|cc|cc|cc|cc|cc|} \hline
 &
  \multicolumn{2}{c|}{$\delta_c=0.5$} &
  \multicolumn{2}{c|}{$\delta_c=1$} &
  \multicolumn{2}{c|}{$\delta_c=2$} &
  \multicolumn{2}{c|}{$\delta_c=4$} &
  \multicolumn{2}{c|}{$\delta_c=6$} &
  \multicolumn{2}{c|}{$\delta_c=8$} \\ \hline
For $x_-$ &
  4 &
  0 &
  4 &
  0 &
  5 &
  0 &
  6 &
  0 &
  10 &
  4 &
  12 &
  5 \\
For $x_+$ &
  4 &
  0 &
  4 &
  0 &
  5 &
  0 &
  6 &
  0 &
  9 &
  1 &
  17 &
  7 \\ \hline
\end{tabular}
\caption{Synthetic dataset: iteration profile of the OptNRQ-nepv (Algorithm~\ref{alg:SCF}) with respect to the tolerance set radius $\delta_c$. Two numbers are displayed for each $\delta_c$, first for the number of iterations and the second for the number of line searches.}
\label{tab:itr_num}
\end{table}

}\end{example}

\begin{example}[Dataset IIIa]
{\rm 
In this example, we investigate a {suitable range of the tolerance set radius} $\delta_c$ values for which the minmax CSP can achieve an improved classification rate compared to the CSP when considering real-world BCI datasets.
We use Dataset IIIa for this investigation.
The trials for each subject are evenly divided between training and testing datasets, as indicated by the provided labels.

Figure~\ref{fig:3a} presents the classification rate results with $\delta_c$ values ranging from 0.1 to 1.5 with a step size of 0.1.
The blue line represents the classification rates of the minmax CSP, solved by theOptNRQ-nepv (Algorithm~\ref{alg:SCF}), for different $\delta_c$ values, while the gray dashed line represents the classification rate of the CSP, serving as a baseline for comparison.

\begin{figure}[tbhp]
    \centering
    \subfloat[]{{\includegraphics[scale=0.14]{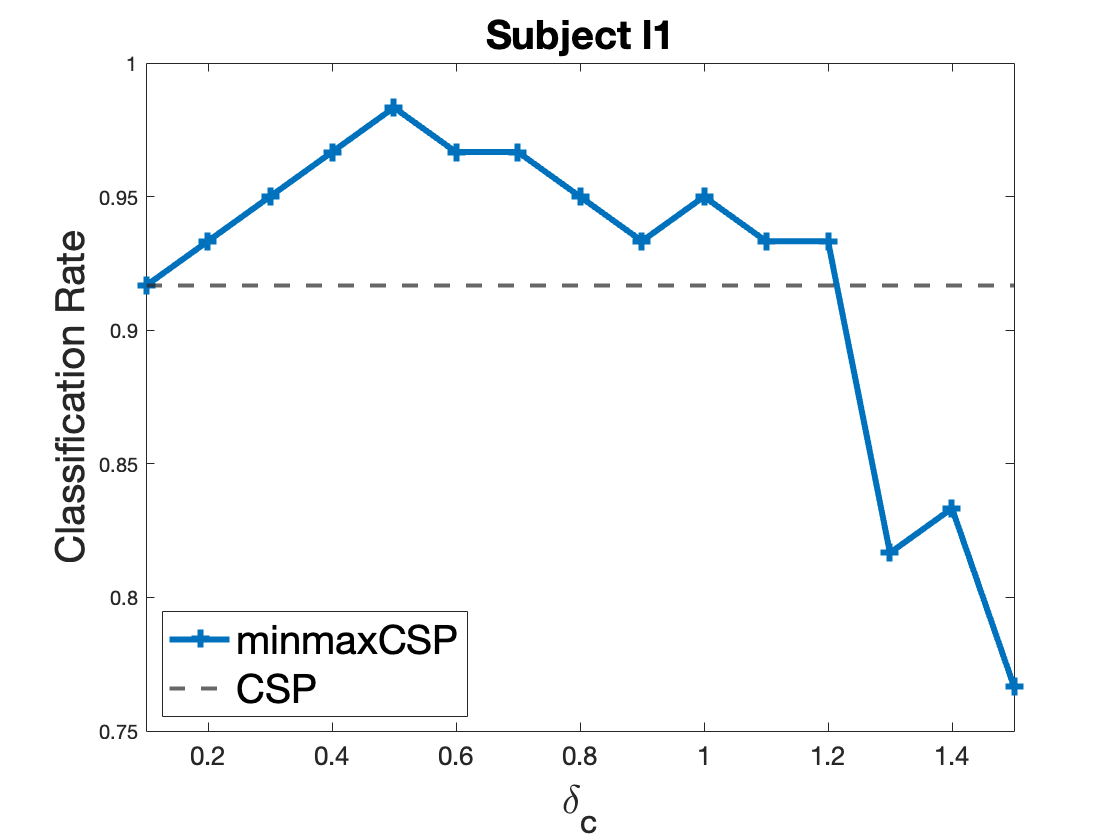} }}%
    \subfloat[]{{\includegraphics[scale=0.14]{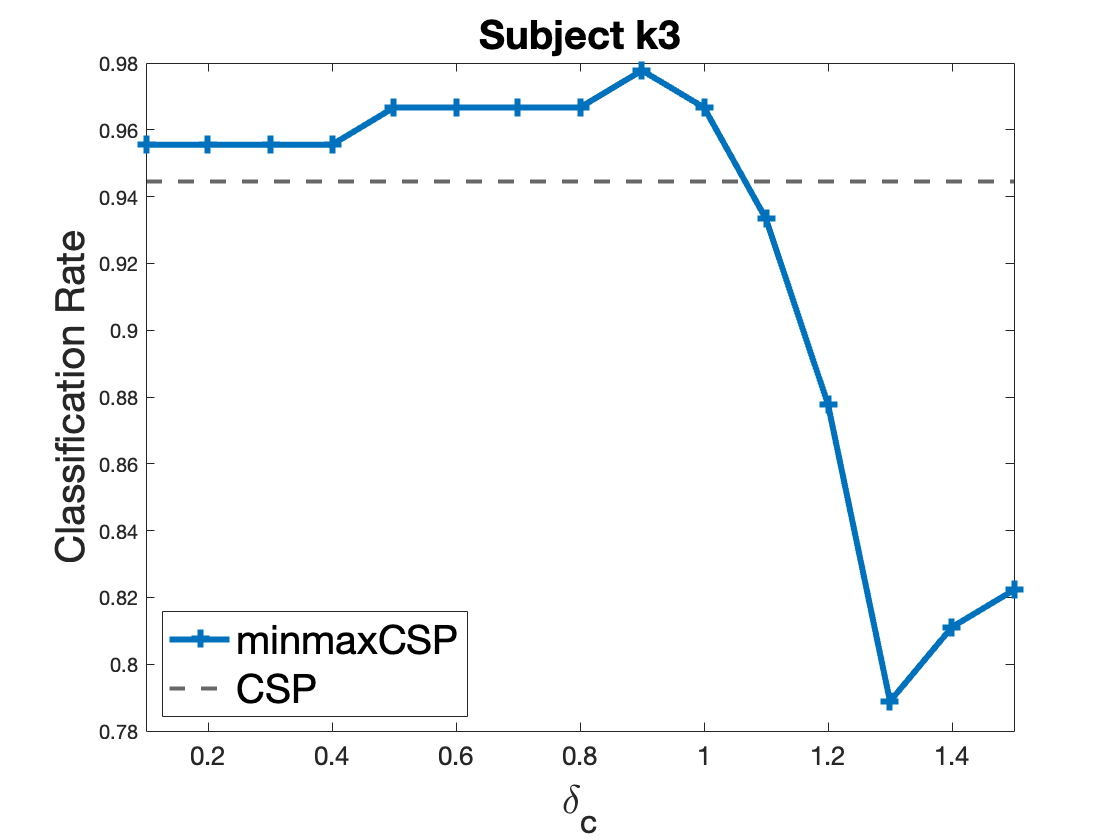} }}%
    \subfloat[]{{\includegraphics[scale=0.14]{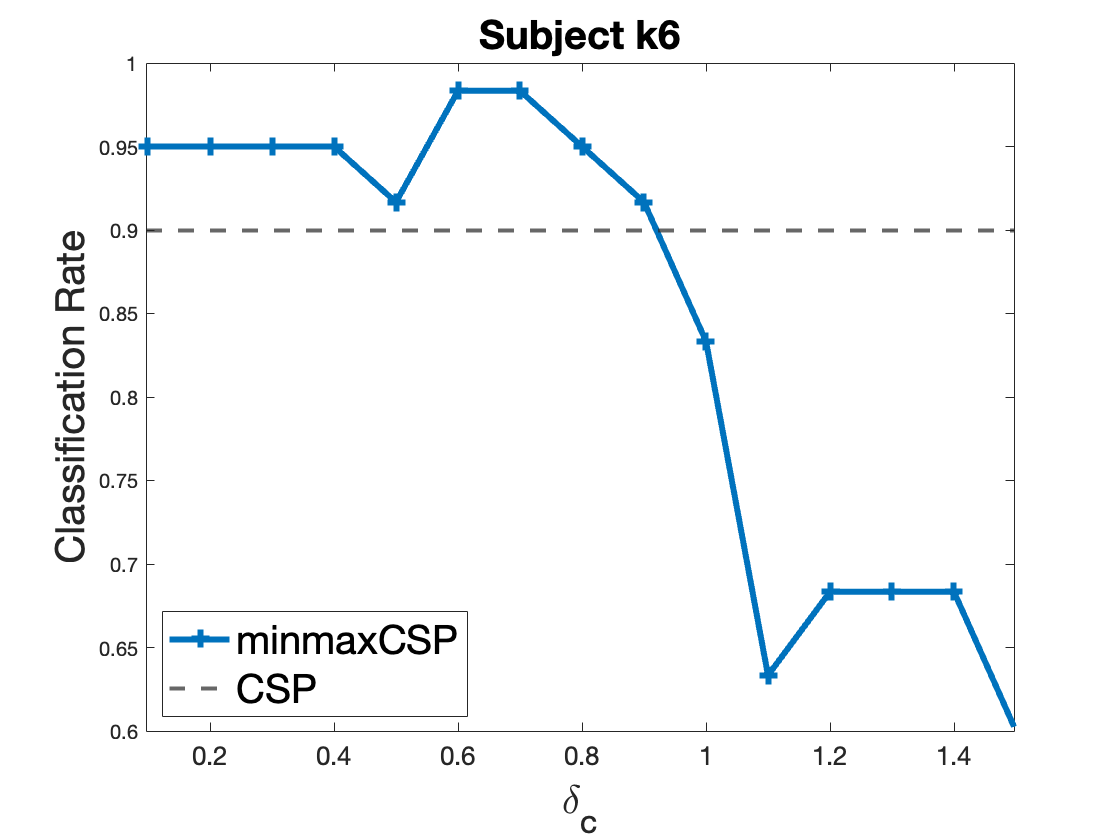} }}%
\caption{Dataset IIIa: classification rate of the OptNRQ-nepv (Algorithm~\ref{alg:SCF}) with respect to the tolerance set radius $\delta_c$ for each subject. The blue line corresponds to the classification rates of the OptNRQ-nepv, while the gray dashed line corresponds to the classification rate of the CSP.}\label{fig:3a}
\end{figure}

Our results indicate that the minmax CSP outperforms the CSP for $\delta_c$ values less than 1. 
However, as the tolerance set radius becomes larger than 1, the performance of the minmax CSP degrades and falls below that of the CSP. 
This decline in performance is attributed to the inclusion of outlying trials within the tolerance set, which results in a poor covariance matrix estimation. 
Based on the results in Figure~\ref{fig:3a}, we restrict the tolerance set 
radii $\delta_c$ to values less than 1 for further 
classification experiments on real-world BCI datasets.
}\end{example}

\begin{example}[Berlin dataset]
{\rm 
For each subject, we use the calibration trials as the training dataset and the feedback trials as the testing dataset. We conduct two types of experiments: (i) comparing the classification performance between the OptNRQ-fp~\eqref{eq:fixed_itr} and the OptNRQ-nepv (Algorithm~\ref{alg:SCF}), and (ii) comparing the classification performance between the CSP and the OptNRQ-nepv.

\begin{enumerate}[label=(\roman*)]

    \item We fix the tolerance set radius $\delta_c$ to be the same for all subjects in both the OptNRQ-fp and the OptNRQ-nepv to ensure a fair comparison. Specifically, we choose $\delta_c=0.4$, a value that yielded the best performance for the majority of subjects in both algorithms. When the OptNRQ-fp fails to converge, we obtain the solutions with the smallest objective value among its iterations. 
The classification rate results are presented in Figure~\ref{fig:Berlin}(a), where each dot represents a subject. The red diagonal line represents equal performance between the OptNRQ-fp and the OptNRQ-nepv, and a dot lying above the red diagonal line indicates an improved classification rate for the OptNRQ-nepv.
We observe that the OptNRQ-nepv outperforms the OptNRQ-fp for half of the subjects. The average classification rate of the subjects is 62.1\% for the OptNRQ-fp, while it is 65.3\% for the OptNRQ-nepv.

    \item For the OptNRQ-nepv, we select the optimal tolerance set radius $\delta_c$ for each subject from a range of 0.1 to 1.0 with a step size of 0.1.
    In this experiment, we observe that 87.5\% of the dots lie above the diagonal line, indicating a significant improvement in classification rates for the OptNRQ-nepv compared to the CSP. 
    The average classification rate of the subjects improved from 64.9\% for the CSP to 70.5\% for the OptNRQ-nepv.
    
\end{enumerate}

\begin{figure}[tbhp]
    \centering
    \subfloat[]{{\includegraphics[scale=0.18]{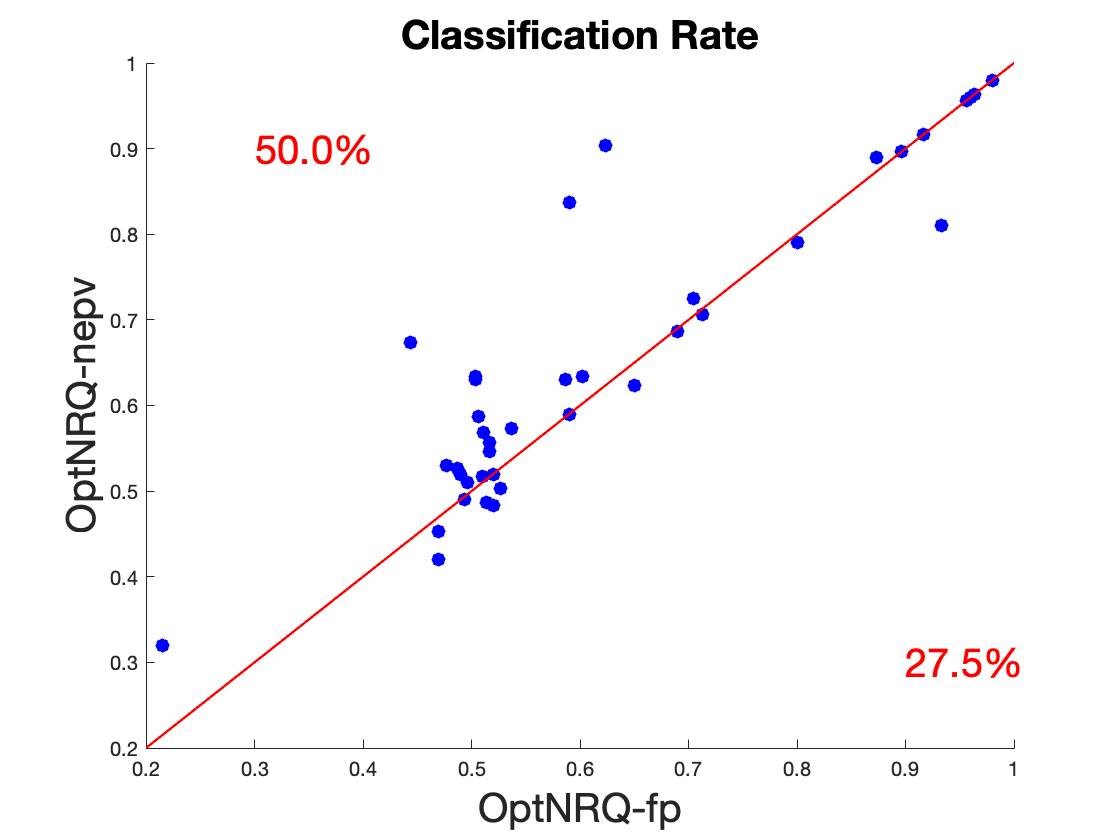} }}%
    \subfloat[]{{\includegraphics[scale=0.18]{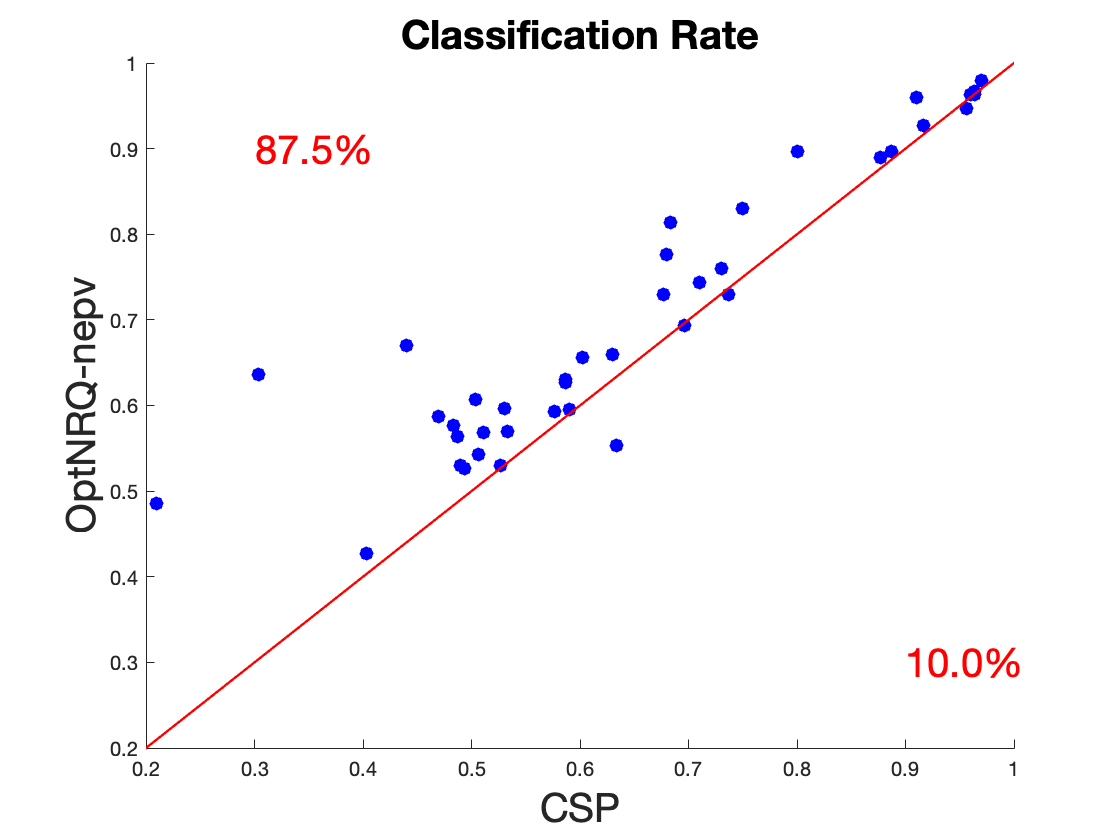} }}%
\caption{Berlin dataset: scatter plots comparing the classification rates. Panel (a) compares the OptNRQ-fp~\eqref{eq:fixed_itr} and the OptNRQ-nepv (Algorithm~\ref{alg:SCF}) and panel (b) compares the CSP and the OptNRQ-nepv. Each dot represents a subject, and if the dot is above the red diagonal line, then the OptNRQ-nepv has a higher classification rate for that subject. The percentages on the upper left and on the lower right indicate the percentage of subjects that perform better for the corresponding algorithm.}\label{fig:Berlin}
\end{figure}

}\end{example}

\begin{example}[Distraction dataset]
{\rm 
In this example, we illustrate the increased robustness of the OptNRQ-nepv (Algorithm~\ref{alg:SCF}) to noises in comparison to the CSP. 
For each subject of the Distraction dataset, we select the trials without any distractions as the training dataset.
The testing dataset is set in six different ways, one for each of the six distraction tasks.
We choose the tolerance set radius $\delta_c$ optimally from a range of 0.1 to 1.0 with a step size of 0.1.


\begin{figure}[tbhp]
    \centering
    \subfloat[]{{\includegraphics[scale=0.12]{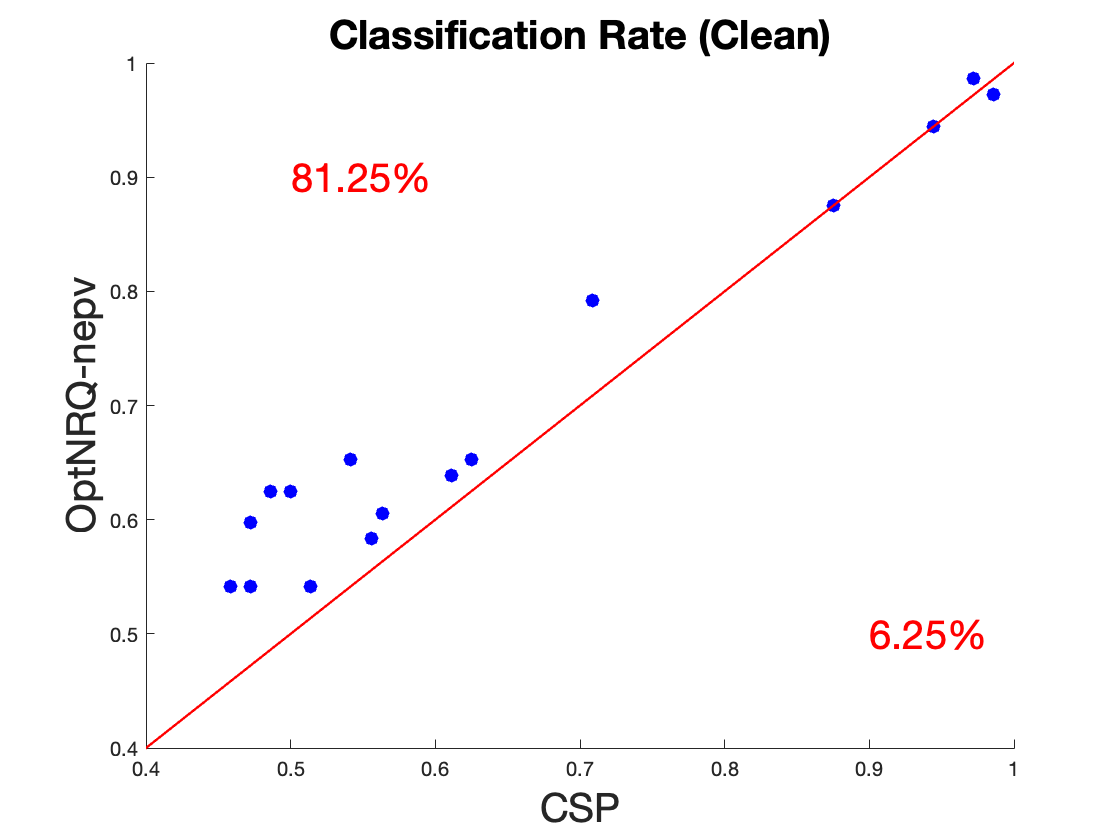} }}%
    \subfloat[]{{\includegraphics[scale=0.12]{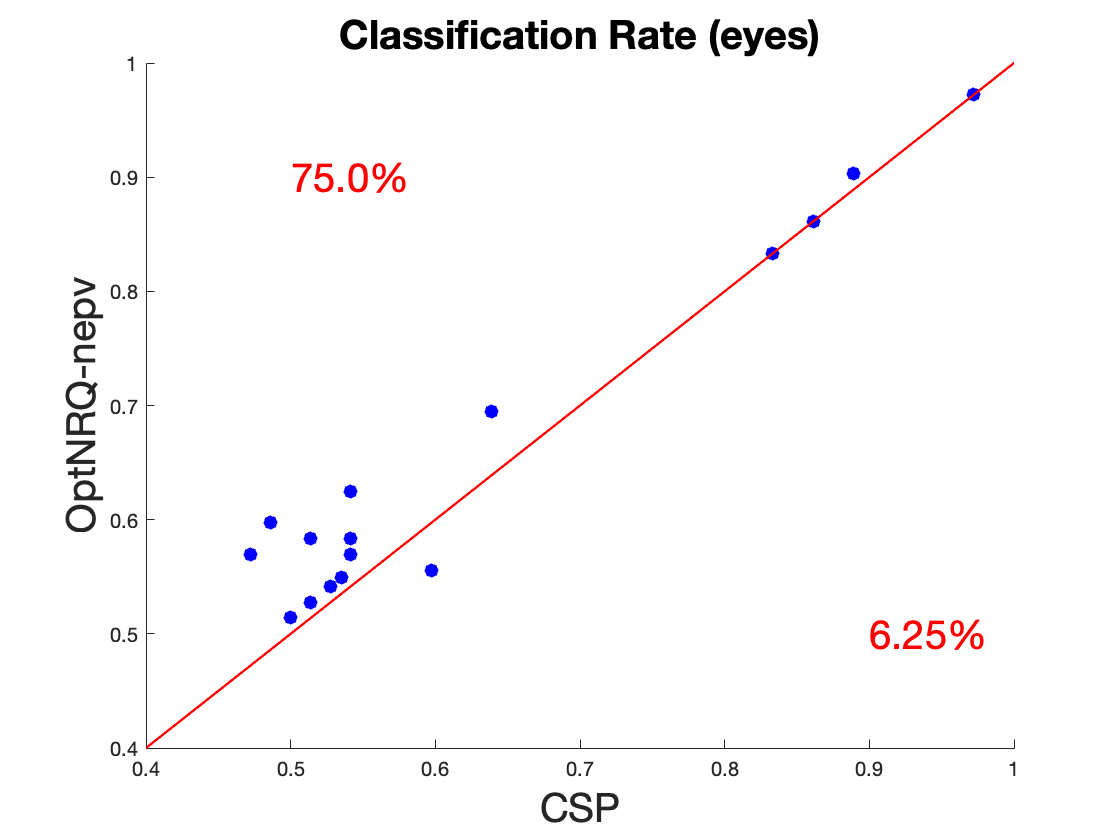} }}%
    \subfloat[]{{\includegraphics[scale=0.12]{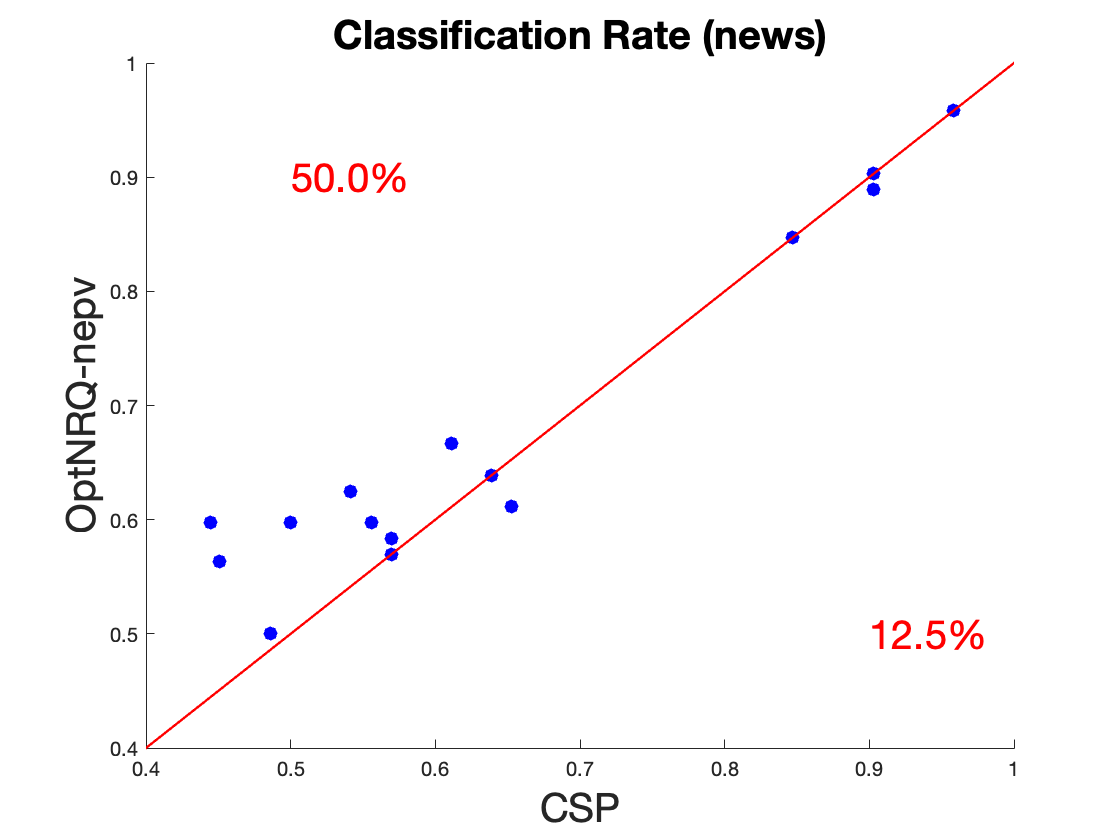} }}%
    \\
    \subfloat[]{{\includegraphics[scale=0.12]{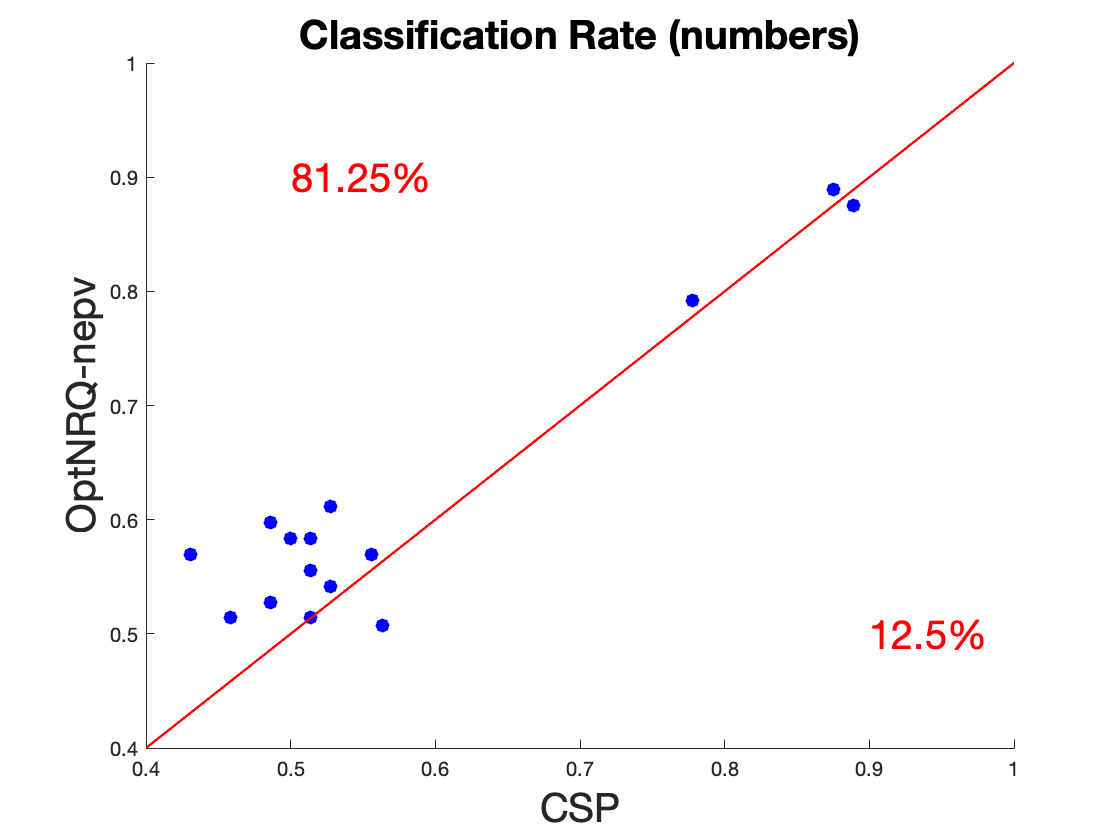} }}%
    \subfloat[]{{\includegraphics[scale=0.12]{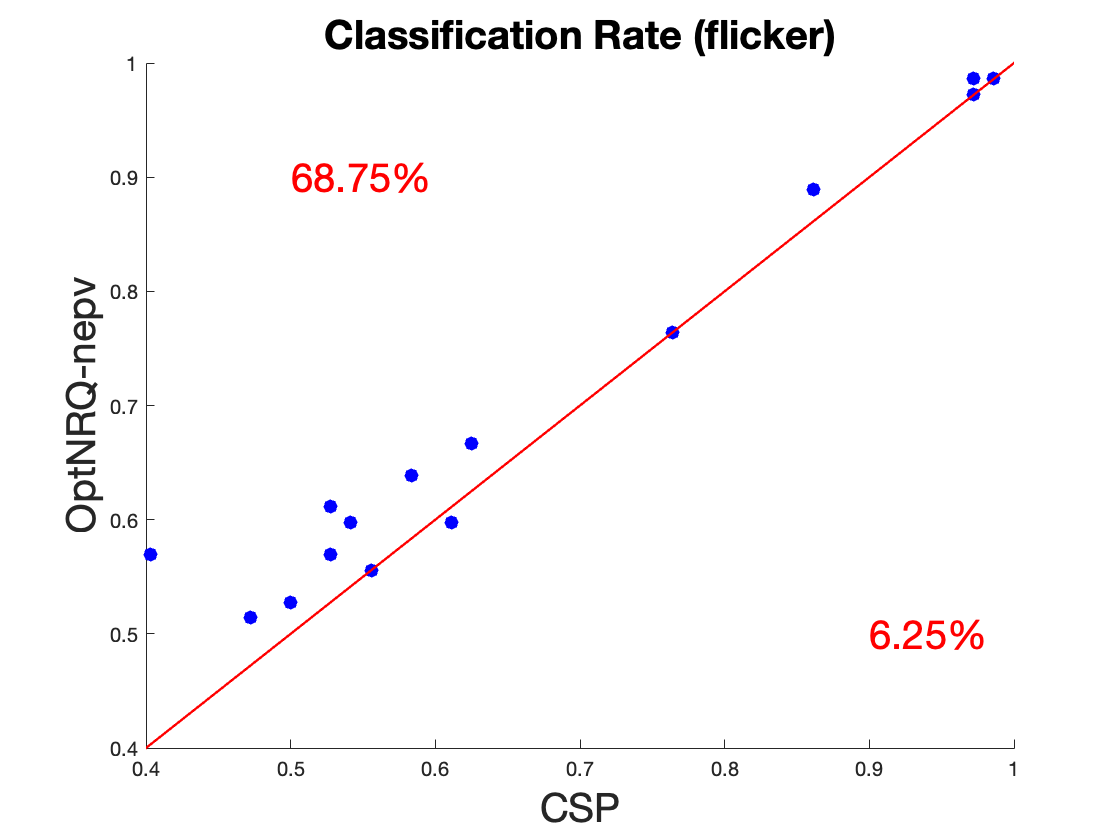} }}%
    \subfloat[]{{\includegraphics[scale=0.12]{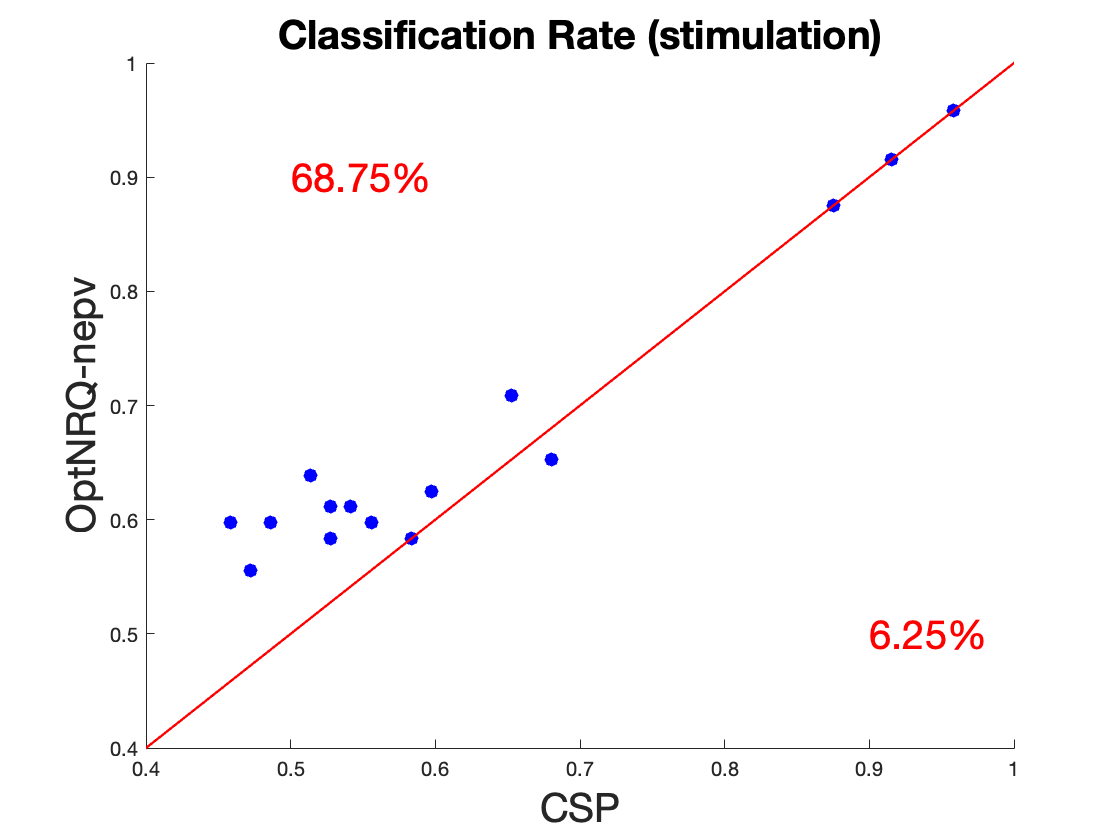} }}%
    \caption{Distraction dataset: scatter plots comparing the classification rates between the CSP and the OptNRQ-nepv (Algorithm~\ref{alg:SCF}). Each plot corresponds to one of the six distraction tasks, where trials of the distraction task is set as the testing dataset.}\label{fig:CSP_Distr}
\end{figure}

The scatter plot in Figure~\ref{fig:CSP_Distr} compares the classification rates of the CSP and the OptNRQ-nepv for each subject.
We observe clear improvements in classification rates for the OptNRQ-nepv compared to the CSP, particularly for the subjects with low BCI control, across all six distraction tasks. 
Overall, the average classification rate of the subjects increased from 62.35\% for the CSP to 66.68\% for the OptNRQ-nepv. 
These results demonstrate the increased robustness of the OptNRQ-nepv to the CSP in the presence of noise in EEG signals.
}\end{example}

\begin{example}[Gwangju dataset]
{\rm 
We follow the cross-validation steps outlined in \cite{cho2017eeg}.
Firstly, the trials for each motor condition are randomly divided into 10 subsets. Next, seven of these subsets are chosen as the training set, and the remaining three subsets are used for testing. 
We compute the classification rates of the CSP and the OptNRQ-nepv (Algorithm~\ref{alg:SCF}), where for the OptNRQ-nepv, we find the optimal tolerance set radius $\delta_c$ for each subject from a range of 0.1 to 1.0 with a step size of 0.1.
For each subject, this procedure is conducted for all 120 possible combinations of training and testing datasets creation. 

\begin{figure}[tbhp]
\centering
\includegraphics[scale=0.18]{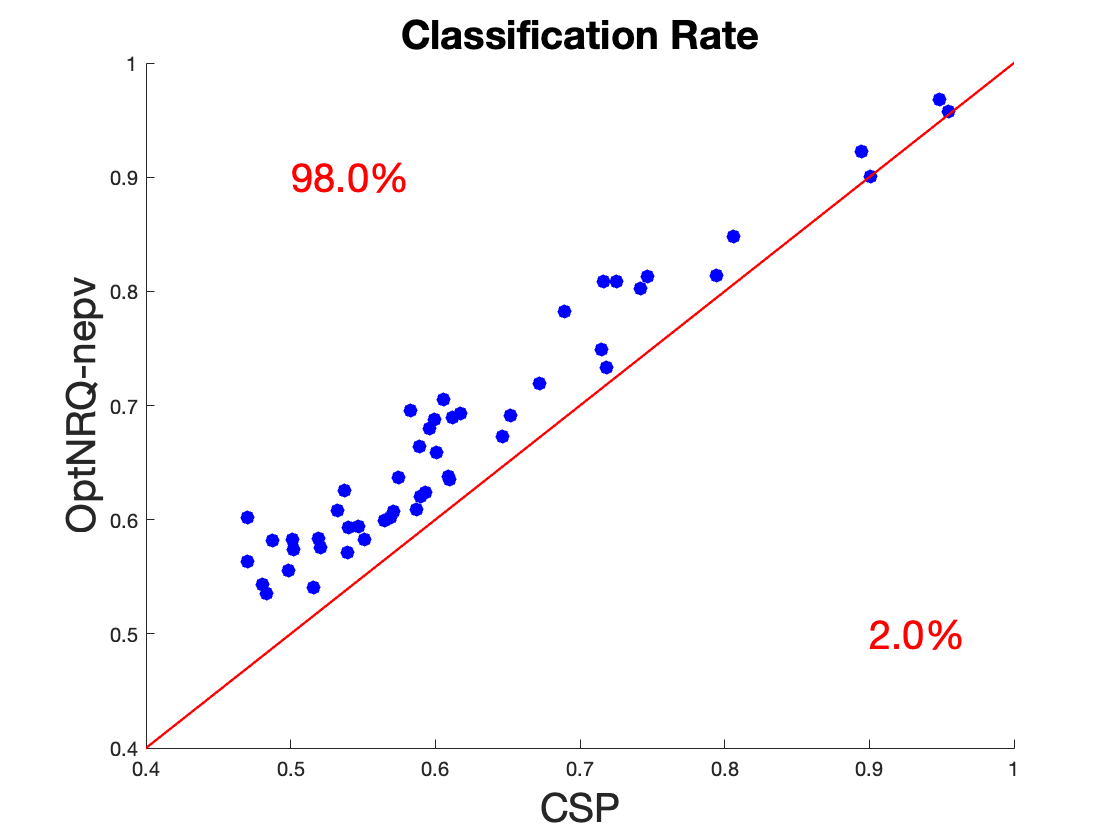}
\caption{Gwangju dataset: scatter plots comparing the classification rates of the CSP and the OptNRQ-nepv (Algorithm~\ref{alg:SCF}).}\label{fig:Gwangju}
\end{figure}

The results for the average classification rates of the 120 combinations are shown in Figure~\ref{fig:Gwangju} for each subject.
With the exception of one subject whose classification rates for the CSP and the OptNRQ-nepv are equal, all other subjects benefit from the OptNRQ-nepv. 
The average classification rate of the subjects improved from 61.2\% for the CSP to 67.0\% for the OptNRQ-nepv.
}\end{example}


\section{Concluding Remarks}\label{sec:con}

We discussed a robust variant of the CSP known as the minmax CSP. 
Unlike the CSP, which utilizes the average covariance matrices that are prone to suffer from poor performance due to EEG signal nonstationarity and artifacts, the minmax CSP considers selecting the covariance matrices from the sets that effectively capture the variability in the signals.
This attributes the minmax CSP to be more robust to nonstationarity and artifacts than the CSP.
The minmax CSP is formulated as a nonlinear Rayleigh quotient optimization.
By utilizing the optimality conditions of this nonlinear Rayleigh quotient optimization, we derived an eigenvector-dependent nonlinear eigenvalue problem (NEPv) with known eigenvalue ordering.
We introduced an algorithm called the OptNRQ-nepv, a self-consistent field iteration with line search, to solve the derived NEPv.
While the existing algorithm for the minmax CSP suffers from convergence issues, we demonstrated that the OptNRQ-nepv converges and displays a local quadratic convergence.
Through a series of numerical experiments conducted on real-world BCI datasets, we showcased the improved classification performance of the OptNRQ-nepv when compared to the CSP and the existing algorithm for the minmax CSP.

We provided empirical results on the impact of the tolerance set radius $\delta_c$ on both the convergence behavior and classification performance of the OptNRQ-nepv. 
While these results provide a suggestion for the range of $\delta_c$ to consider when working with real-world EEG datasets, we acknowledge that determining the optimal $\delta_c$ remains an open question.
Our future studies will explore a systematic approach for calculating the optimal value of $\delta_c$.

We also acknowledge the limitations of the considered data-driven tolerance sets in the minmax CSP.
These tolerance sets are computed using Euclidean-based PCA on the data covariance matrices, potentially failing to precisely capture their variability within the non-Euclidean matrix manifold.
To address this, we plan to investigate the application of different tolerance sets, such as non-convex sets, which might more accurately capture the variability of the data covariance matrices.

It is common in the CSP to compute multiple spatial filters (i.e., $k>1$ many smallest eigenvectors of \eqref{eq:gen_eig}).
Using multiple spatial filters often leads to improved classification rates when compared with that of using a single spatial filter.
We can aim to compute multiple spatial filters with the minmax CSP as well.
This amounts to extending the Rayleigh quotient formulation of the minmax CSP~\eqref{eq:cspminmax-2} to arrive at a minmax trace ratio optimization.
Similar to the derivation of the nonlinear Rayleigh quotient optimization~\eqref{eq:cspnrq} from the minmax CSP~\eqref{eq:cspminmax-2} shown in this paper, it can be shown that a nonlinear trace ratio optimization can be derived from the extended minmax trace ratio optimization.
However, the methodology for solving the nonlinear trace ratio optimization remains an active research problem.
Our future plans involve a comprehensive review of related research on nonlinear trace ratio optimization, such as Wasserstein discriminant analysis \cite{flamary2018wasserstein, liu2020ratio, roh2022bi}, to explore various approaches for addressing this problem.







\begin{appendices}

\section{Frobenius-norm Tolerance Sets}\label{app:Frob}

In this approach, the tolerance sets $\mathcal{S}_c$ are constructed as the ellipsoids of covariance matrices centered around the average covariance matrices $\widehat{\Sigma}_c$ \cite{kawanabe2009robust,kawanabe2014robust}. 
They are
\begin{equation}\label{eq:gen_tol}
    \mathcal{S}_c=\{\Sigma_c\mid \Sigma_c\succ0,\|\Sigma_c-\widehat{\Sigma}_c\|_P\leq\delta_c\},
\end{equation}
where $\|\cdot\|_P$ is the {weighted Frobenius norm}
\begin{equation}\label{eq:wF_norm}
    \|X\|_P := \|P^{1/2}XP^{1/2}\|_F=\sqrt{\tr(PX^TPX)}
\end{equation}
defined for $X\in\R^{n\times n}$ with $P\in\R^{n\times n}$ as a symmetric positive definite matrix. 
$\|\cdot\|_F$ denotes the Frobenius norm. 
The parameter $\delta_c$ denotes the radius of the ellipsoids. 

In the following Lemma, we show that the worst-case covariance matrices in the tolerance sets~\eqref{eq:gen_tol} can be determined explicitly for any choice of positive definite matrix $P$ in the norm~\eqref{eq:wF_norm}.

\begin{lemma}[{\cite[Lemma 1]{kawanabe2014robust}}]\label{lem:P}
Let the tolerance sets $\mathcal{S}_c$ and $\mathcal{S}_{\bar{c}}$ be defined as~\eqref{eq:gen_tol} with the norm~\eqref{eq:wF_norm} at $P=P_c^{-1}$ and $P=P_{\bar{c}}^{-1}$, respectively. 
Then, the minmax optimization~\eqref{eq:cspminmax-} becomes
\begin{align}\label{eq:csp_P_norm-}
    \min_{x\neq0}
    \frac{x^T(\widehat{\Sigma}_c+\delta_cP_c)x}{x^T(\widehat{\Sigma}_c+\widehat{\Sigma}_{\bar{c}}+\delta_cP_c-\delta_{\bar{c}}P_{\bar{c}})x}.
\end{align}
\end{lemma}

\begin{proof}
The optimal $\Sigma_c$ and $\Sigma_{\bar{c}}$ that maximize the inner optimization of the Rayleigh quotient~\eqref{eq:cspminmax-} can be found from two separate optimizations. 
Since $\Sigma_c$ and $\Sigma_{\bar{c}}$ are positive definite, we have
\[\max_{\substack{\Sigma_c\in\mathcal{S}_c\\ \Sigma_{\bar{c}}\in\mathcal{S}_{\bar{c}}}}\frac{x^T\Sigma_cx}{x^T(\Sigma_c+\Sigma_{\bar{c}})x}=\frac{\max\limits_{\Sigma_c\in\mathcal{S}_c}x^T\Sigma_cx}{\max\limits_{\Sigma_c\in\mathcal{S}_c}x^T\Sigma_cx + \min\limits_{\Sigma_{\bar{c}}\in\mathcal{S}_{\bar{c}}}x^T\Sigma_{\bar{c}}x}.\]
So, we need to show the equality for two optimizations
\begin{equation}\label{eq:first_opt}
    \argmax\limits_{\Sigma_c\in\mathcal{S}_c}x^T\Sigma_cx = \widehat{\Sigma}_c+\delta_cP_c
\end{equation}
and
\begin{equation}\label{eq:second_opt}
    \argmin\limits_{\Sigma_{\bar{c}}\in\mathcal{S}_{\bar{c}}}x^T\Sigma_{\bar{c}}x = \widehat{\Sigma}_{\bar{c}}-\delta_{\bar{c}}P_{\bar{c}}.
\end{equation}
These two optimizations are equivalent to
\begin{align*}
    \argmax\limits_{\Sigma_c\in\mathcal{S}_c}x^T\Sigma_cx &= \argmax\limits_{\Sigma_c\in\mathcal{S}_c}x^T(\Sigma_c - \widehat{\Sigma}_c)x
\end{align*}
and
\begin{align*}
    \argmin\limits_{\Sigma_{\bar{c}}\in\mathcal{S}_{\bar{c}}}x^T\Sigma_{\bar{c}}x &= \argmin\limits_{\Sigma_{\bar{c}}\in\mathcal{S}_{\bar{c}}}x^T(\Sigma_{\bar{c}} - \widehat{\Sigma}_{\bar{c}})x,
\end{align*}
respectively.

We prove the equality of~\eqref{eq:first_opt} by providing an upper bound of $|x^T(\Sigma_c-\widehat{\Sigma}_c)x|$ and by determining $\Sigma_c$ that achieves such upper bound.
To do so, we make use of the inner product that induces the weighted Frobenius norm~\eqref{eq:wF_norm}.
This inner product is defined as
\begin{equation}\label{eq:inn_wF}
    \langle X,Y \rangle_P := \tr(PX^TPY).
\end{equation}
The Cauchy-Schwarz inequality corresponding to the inner product~\eqref{eq:inn_wF} is
\begin{equation}\label{eq:CS}
    |\langle X,Y \rangle_P| \leq \|X\|_P \|Y\|_P = \sqrt{\tr(PX^TPX)} \cdot \sqrt{\tr(PY^TPY)}.
\end{equation}
We have 
\begin{align}
    |x^T(\Sigma_c-\widehat{\Sigma}_c)x|&=|\tr((\Sigma_c-\widehat{\Sigma}_c)xx^T)|=|\tr(P_cP_c^{-1}(\Sigma_c-\widehat{\Sigma}_c)P_c^{-1}P_cxx^T)| \nonumber \\
    &= |\tr(P_c^{-1}(\Sigma_c-\widehat{\Sigma}_c)^TP_c^{-1}(P_cxx^TP_c))| = |\langle (\Sigma_c-\widehat{\Sigma}_c), (P_cxx^TP_c) \rangle_{P_c^{-1}}| \nonumber \\
    &\leq \sqrt{\tr(P_c^{-1}(\Sigma_c-\widehat{\Sigma}_c)^TP_c^{-1}(\Sigma_c-\widehat{\Sigma}_c))} \cdot \sqrt{\tr(P_c^{-1}(P_cxx^TP_c)^TP_c^{-1}(P_cxx^TP_c))} \label{eq:lem1_ineq1} \\
    &\leq \delta_c \cdot \sqrt{\tr((x^TP_cx)^2)} =  \delta_c \cdot x^TP_cx \label{eq:lem1_ineq2}
\end{align}
where the inequality~\eqref{eq:lem1_ineq1} follows from the Cauchy-Schwarz inequality~\eqref{eq:CS}, and the inequality~\eqref{eq:lem1_ineq2} follows from the constraint $\|\Sigma_c-\widehat{\Sigma}_c\|_{P_c^{-1}} \leq\delta_c$ of the tolerance set $\mathcal{S}_c$~\eqref{eq:gen_tol}.
With $\Sigma_c=\widehat{\Sigma}_c + \delta_cP_c$, we have $|x^T(\Sigma_c-\widehat{\Sigma}_c)x| = \delta_c \cdot x^TP_cx$. 
Therefore, for the optimization~\eqref{eq:first_opt} the maximization is reached by $\Sigma_c=\widehat{\Sigma}_c+\delta_cP_c$. 

By analogous steps, and with an assumption that $\widehat{\Sigma}_{\bar{c}}-\delta_{\bar{c}}P_{\bar{c}}$ is positive definite, the minimization of~\eqref{eq:second_opt} is reached by $\Sigma_{\bar{c}}=\widehat{\Sigma}_{\bar{c}}-\delta_{\bar{c}}P_{\bar{c}}$.

\end{proof}

\section{Vector and Matrix Calculus}\label{app:der}

Let $x$ be a vector of size $n$:
$$x=\begin{bmatrix}
    x_1\\
    x_2\\
    \vdots\\
    x_n
\end{bmatrix}\in\R^n.
$$
Denoting $\nabla$ as a gradient operator with respect to $x$, the gradient of a differentiable function $a(x):\R^n\to\R$ with respect to $x$ is defined as
\begin{equation}\label{eq:a_der}
    \nabla a(x):=\begin{bmatrix}
    \frac{\partial a(x)}{\partial x_1}\\
    \frac{\partial a(x)}{\partial x_2}\\
    \vdots\\
    \frac{\partial a(x)}{\partial x_n}\\
\end{bmatrix}\in\R^n.
\end{equation}
Similarly, the gradient of a differentiable vector-valued function $y(x):\R^n\to\R^m$ with respect to $x$ is defined as
\begin{align}\label{eq:y_der}
    \nabla y(x)&:=\begin{bmatrix}
    \frac{\partial y_1(x)}{\partial x_1} & \frac{\partial y_2(x)}{\partial x_1} & \cdots & \frac{\partial y_m(x)}{\partial x_1}\\
    \frac{\partial y_1(x)}{\partial x_2} & \frac{\partial y_2(x)}{\partial x_2} & \cdots & \frac{\partial y_m(x)}{\partial x_2}\\
    \vdots & \vdots & \mbox{ } & \vdots \\
    \frac{\partial y_1(x)}{\partial x_n} & \frac{\partial y_2(x)}{\partial x_n} & \cdots & \frac{\partial y_m(x)}{\partial x_n}\\
\end{bmatrix} \nonumber \\ 
    &\mbox{ }= \begin{bmatrix}
        \nabla y_1(x) & \nabla y_2(x) & \cdots & \nabla y_m(x)
    \end{bmatrix} \in\R^{n\times m},
\end{align}
where $\nabla y_i(x)\in\R^n$ is the gradient of the $i^{\mbox{th}}$ component of $y(x)$.

In the following two lemmas, we show vector and matrix gradients that are utilized in Section~\ref{sec:NEPv}. 
The first lemma consists of basic properties of gradients whose proofs can be found in a textbook such as \cite[Appendix D]{dattorro2010convex}. 
For the second lemma, we provide the proofs.

\begin{lemma}\label{lem:der1}
    Let $A\in\R^{n\times n}$ be a symmetric matrix. Then, the following gradients hold:
    \begin{enumerate}[label=(\alph*)]
        \item $\nabla(Ax)=A$.
        \item $\nabla(x^TAx)=2Ax$.
    \end{enumerate}
\end{lemma}

\begin{lemma}\label{lem:der2}
    Let $a(x):\R^n\to\R$, $y(x):\R^n\to\R^m$, and $B(x):\R^n\to\R^{n\times m}$ be differentiable real-valued function, vector-valued function, and matrix-valued function of $x$, respectively. Then, the following gradients hold:
    \begin{enumerate}[label=(\alph*)]
        \item Given a diagonal matrix $W\in\R^{m\times m}$, $\nabla(y(x)^TWy(x))=2\nabla y(x)Wy(x)\in\R^n$.
        \item $\nabla(B(x)y(x))=\sum_{i=1}^my_i(x)\nabla b_i(x) + \nabla y(x)B(x)^T\in\R^{n\times n}$, where $b_i(x)$, for $i=1,2,\ldots,m$, represent the column vectors of $B(x)$.
        \item $\nabla(a(x)y(x))=\nabla a(x)y(x)^T+a(x)\nabla y(x)\in\R^{n\times m}$.
    \end{enumerate}
\end{lemma}
\begin{proof}
(a) We have $y(x)^TWy(x)=\sum_{i=1}^mw_i(y_i(x))^2$ where $w_i$ is the $i^{\mbox{th}}$ diagonal element of $W$. Then, by the chain rule,
        $$\nabla(y(x)^TWy(x))=\sum_{i=1}^m2w_iy_i(x)\nabla y_i(x),$$
        which is precisely $2\nabla y(x)Wy(x)$ by the definition~\eqref{eq:y_der} of $\nabla y(x)$.

(b) We have 
        $$B(x)y(x)=\begin{bmatrix}
            b_{11}(x) & \cdots & b_{1m}(x) \\
            \vdots & \mbox{ } & \vdots \\
            b_{n1}(x) & \cdots & b_{nm}(x)
        \end{bmatrix}
        \begin{bmatrix}
            y_1(x)\\
            \vdots\\
            y_m(x)
        \end{bmatrix}=
        \begin{bmatrix}
            \sum_{j=1}^mb_{1j}(x)y_j(x) \\
            \vdots \\
            \sum_{j=1}^mb_{nj}(x)y_j(x)
        \end{bmatrix}\in\R^n,$$
        and the $i^{\mbox{th}}$ column of the gradient $\nabla(B(x)y(x))\in\R^{n\times n}$ is the gradient of $B(x)y(x)$'s $i^{\mbox{th}}$ element:
        $$\nabla(\sum_{j=1}^mb_{ij}(x)y_j(x))=\sum_{j=1}^m y_j(x)\nabla b_{ij}(x) + \sum_{j=1}^mb_{ij}(x)\nabla y_j(x).$$
        $\nabla(B(x)y(x))$ can then be written as a sum of two matrices
        $$\nabla(B(x)y(x))=\mathcal{B}(x) + \mathcal{Y}(x),$$
        where 
        $$\mathcal{B}(x) := \begin{bmatrix}
            \sum_{j=1}^m y_j(x)\nabla b_{1j}(x) & \sum_{j=1}^m y_j(x)\nabla b_{2j}(x) & \cdots & \sum_{j=1}^m y_j(x)\nabla b_{nj}(x)
        \end{bmatrix}\in\R^{n\times n}$$
        and
        $$\mathcal{Y}(x) := \begin{bmatrix}
            \sum_{j=1}^mb_{1j}(x)\nabla y_j(x) & \sum_{j=1}^mb_{2j}(x)\nabla y_j(x) & \cdots & \sum_{j=1}^mb_{nj}(x)\nabla y_j(x)
        \end{bmatrix}\in\R^{n\times n}.$$
        By decomposing, $\mathcal{B}(x)$ is equal to the sum of $m$ many $n\times n$ matrices:
        \begin{align*}
            \mathcal{B}(x) &= \begin{bmatrix}
            y_1(x)\nabla b_{11}(x) & y_1(x)\nabla b_{21}(x) & \cdots & y_1(x)\nabla b_{n1}(x)
        \end{bmatrix} \\
        &\mbox{    } + \begin{bmatrix}
            y_2(x)\nabla b_{12}(x) & y_2(x)\nabla b_{22}(x) & \cdots & y_2(x)\nabla b_{n2}(x)
        \end{bmatrix} \\
        &\mbox{    } + \cdots + \begin{bmatrix}
            y_m(x)\nabla b_{1m}(x) & y_m(x)\nabla b_{2m}(x) & \cdots & y_m(x)\nabla b_{nm}(x)
        \end{bmatrix} \\
        & = y_1(x)\begin{bmatrix}
            \nabla b_{11}(x) & \nabla b_{21}(x) & \cdots & \nabla b_{n1}(x)
        \end{bmatrix} \\
        &\mbox{    } + y_2(x)\begin{bmatrix}
            \nabla b_{12}(x) & \nabla b_{22}(x) & \cdots & \nabla b_{n2}(x)
        \end{bmatrix} \\
        &\mbox{    } + \cdots +y_m(x)\begin{bmatrix}
            \nabla b_{1m}(x) & \nabla b_{2m}(x) & \cdots & \nabla b_{nm}(x)
        \end{bmatrix} \\
        &= \sum_{i=1}^my_i(x)\nabla b_i(x)
        \end{align*}
        where $\nabla b_i(x)$ is the gradient of the $i^{\mbox{th}}$ column of $B(x)$.

        For $\mathcal{Y}(x)$, note that 
        \[\nabla y(x)B(x)^T = \sum_{j=1}^m\nabla y_j(x)b_j(x)^T.\]
        This indicates that the $i^{\mbox{th}}$ column of $\nabla y(x)B(x)^T$ corresponds to $\sum_{j=1}^mb_{ij}(x)\nabla y_j(x)$, which is the $i^{\mbox{th}}$ column of $\mathcal{Y}(x)$. Hence, we have $\mathcal{Y}(x)=\nabla y(x)B(x)^T$.

(c) We have 
        $$a(x)y(x)=\begin{bmatrix}
            a(x)y_1(x) \\
            a(x)y_2(x) \\
            \vdots \\
            a(x)y_m(x)
        \end{bmatrix}$$
        so that
        \begin{align*}
            \nabla(a(x)y(x))&=\begin{bmatrix}
            \nabla(a(x)y_1(x)) & \cdots & \nabla(a(x)y_m(x)) 
        \end{bmatrix} \\
        &=\begin{bmatrix}
            y_1(x)\nabla a(x)+a(x)\nabla y_1(x) & \cdots & y_m(x)\nabla a(x)+a(x)\nabla y_m(x)
        \end{bmatrix} \\
        &=\begin{bmatrix}
            y_1(x)\nabla a(x) & \cdots & y_m(x)\nabla a(x)
        \end{bmatrix} +
        \begin{bmatrix}
            a(x)\nabla y_1(x) & \cdots & a(x)\nabla y_m(x)
        \end{bmatrix} \\
        &=\nabla a(x)y(x)^T+a(x)\nabla y(x).
        \end{align*}
\end{proof}

\end{appendices}


\bibliographystyle{plain}
\bibliography{refs}

\end{document}